\documentclass[10pt, letterpaper]{article}

\usepackage{graphicx}
\usepackage{epsfig} 
\usepackage{amsmath} 
\usepackage{amssymb}  
\usepackage{mathrsfs}
\usepackage{url}
\usepackage{tabularx,ragged2e,booktabs,caption}
\usepackage{bm}
\usepackage{bbm}
\usepackage{color}
\usepackage{theorem}
\usepackage{caption}
\usepackage{subcaption}
\usepackage{supertabular}
\usepackage{hyperref}
\usepackage{float}
\usepackage{longtable}
\usepackage{booktabs}
\usepackage{array}
\usepackage{url}

\allowdisplaybreaks

\newcommand{\be}[1]{\begin{equation}\label{#1}}
\newcommand{\benon}{\begin{equation*}}  
\newcommand{\bemuln}[1]{\begin{multline}\label{#1}}
\newcommand{\bemul}{\begin{multline*}}
\newcommand{\bee}{\begin{eqnarray*}}
\newcommand{\eee}{\end{eqnarray*}}
\newcommand{\been}[1]{\begin{eqnarray}\label{#1}}
\newcommand{\eeen}{\end{eqnarray}}
\newcommand{\began}[1]{\begin{gather}\label{#1}}
\newcommand{\bega}{\begin{gather*}}
\newcommand{\bealn}[1]{\begin{align}\label{#1}}
\newcommand{\beal}{\begin{align*}}
\newcommand{\bealatn}[2]{\begin{alignat}{#1}\label{#2}}
\newcommand{\bealat}{\begin{alignat*}}
\newcommand{\bexalatn}[1]{\begin{xalignat}\label{#1}}
\newcommand{\bexalat}{\begin{xalignat*}}

\newcommand{\qed}{\newline \mbox{ } \hfill 
            \rule[-1pt]{2.5mm}{2.5mm}\par\vskip 10pt }

\newcommand{\pf}{\vskip 5 pt \noindent {\bf{Proof: }}}

\newcommand{\mb}{\mathbf}

{\theoremstyle{plain} \newtheorem{thm}{Theorem}[section]

\newtheorem{prop}[thm]{Proposition}
}
{\theoremstyle{break} \theorembodyfont{\it}
\newtheorem{defi}{Definition}
\newtheorem{ass}{Assumption} 
\newtheorem{rmk}{Remark} 
 }

\def\bb{{\mathbf b}}

\def\bd{{\mathbf d}}

\def\bg{{\mathbf g}}
\def\bh{{\mathbf h}}

\def\bm{{\mathbf m}}

\def\bt{{\mathbf t}}

\def\bw{{\mathbf w}}
\def\bx{{\mathbf x}}  
\def\by{{\mathbf y}}

\def\bA{{\mathbf A}}

\def\bM{{\mathbf M}}
\def\bN{{\mathbf N}}

\def\bP{{\mathbf P}}
\def\bQ{{\mathbf Q}}

\def\bS{{\mathbf S}}

\def\texitem#1{\par\smallskip\noindent\hangindent 25pt
               \hbox to 25pt {\hss #1 ~}\ignorespaces}

\newcommand{\bzero}{{\mathbf{0}}}

\newcommand{\scrA}{\mathcal{A}}

\newcommand{\scrD}{\mathcal{D}}

\newcommand{\scrF}{\mathcal{F}}

\newcommand{\scrH}{\mathcal{H}}
\newcommand{\scrI}{\mathcal{I}}

\newcommand{\scrK}{\mathcal{K}}

\newcommand{\scrR}{\mathcal{R}}
\newcommand{\scrS}{\mathcal{S}}

\newcommand{\scrV}{\mathcal{V}}
\newcommand{\scrW}{\mathcal{W}}

\newcommand{\bbeta}{\boldsymbol{\beta}}

\newcommand{\bepsilon}{\boldsymbol{\epsilon}}

\newcommand{\bxi}{\boldsymbol{\xi}}

\usepackage{algorithm}
\usepackage{algorithmicx,algpseudocode}

\DeclareMathOperator{\VI}{\text{VI}}

\makeatletter
\newcommand*{\defeq}{\stackrel{\text{def}}{=}}
\makeatother

\makeatletter
\newcommand{\rmnum}[1]{\romannumeral #1}
\newcommand{\Rmnum}[1]{\expandafter\@slowromancap\romannumeral #1@}
\makeatother

\begin{document}

\title{The Price of Anarchy in Transportation Networks: Data-Driven
  Evaluation and Reduction Strategies \thanks{Research partially
    supported by the NSF under grants CNS-1645681, ECCS-1509084,
    CCF-1527292, IIS-1237022, and IIP-1430145, by the DOE under grant
    DE-AR0000796, by the AFOSR under grant FA9550-15-1-0471, by the ARO
    under grant W911NF-12-1-0390, by Bosch, and by The MathWorks.}}

\author{Jing Zhang, Sepideh Pourazarm, \thanks{$^\dag$ Division of
    Systems Eng., Boston University, Boston, MA 02446, email: {\tt
      \{jzh, sepid\}@bu.edu}} Christos G. Cassandras,\\
and Ioannis Ch. Paschalidis \thanks{$^\ddag$ Dept. of Electrical and
  Computer Eng. and Division of Systems Eng.,  
Boston University, 8 St. Mary's St., Boston, MA 02215, email: {\tt
 \{cgc, yannisp\}@bu.edu}, url: {\tt \url{http://sites.bu.edu/paschalidis}}.}
}

\maketitle
\begin{abstract}
  Among the many functions a Smart City must support, transportation
  dominates in terms of resource consumption, strain on the environment,
  and frustration of its citizens. We study transportation networks
  under two different routing policies, the commonly assumed selfish
  user-centric routing policy and a socially-optimal system-centric
  one. We consider a performance metric of efficiency -- the Price of
  Anarchy (PoA) -- defined as the ratio of the total travel latency cost
  under selfish routing over the corresponding quantity under
  socially-optimal routing. We develop a data-driven approach to
  estimate the PoA, which we subsequently use to conduct a case study
  using extensive actual traffic data from the Eastern Massachusetts
  road network. To estimate the PoA, our approach learns from data a
  complete model of the transportation network, including
  origin-destination demand and user preferences. We leverage this model
  to propose possible strategies to reduce the PoA and increase
  efficiency.
\end{abstract}


\section{Introduction} \label{sec:intro}

As of 2014, 54\% of the earth's population resides in urban areas, a
percentage expected to reach 66\% by 2050. This increase would amount to
2.5 billion people added to urban populations~\cite{un-desa-14}. At the
same time, there are now 28 mega-cities (with population $\geq 10$
million) worldwide, accounting for 22\% of the world's urban
dwellers, and projections indicate more than 41 mega-cities by 2030. It
stands to reason that the management and sustainability of urban areas
has become one of the most critical challenges our societies face today,
leading to a quest for ``smart'' cities.

Among the many functions a city supports, transportation dominates in
terms of resource consumption, strain on the environment, and
frustration of its citizens. Commuter delays have risen by 260\% over
the past 25 years and 28\% of U.S. primary energy is now used in
transportation~\cite{schrank2011tti}. It is estimated that the
cumulative cost of traffic congestion by 2030 will reach \$2.8
trillion~\cite{inrix-2015} -- equal roughly to the U.S. annual tax
revenue. This estimate accounts for direct costs to drivers (time, fuel)
and indirect costs resulting from businesses passing these same costs on
to consumers, but it does not include the equally alarming environmental
impact due to a large proportion of toxic air pollutants attributed to
mobile sources. At the individual citizen level, traffic congestion led
to \$1,740 in average costs per driver during 2014. If unchecked, this
number is expected to grow by more than 60\%, to \$2,900 annually, by
2030~\cite{inrix-2015}.

A transportation network is a system with non-cooperative agents
(drivers) in which each agent seeks to minimize her own individual cost
by choosing the best route (resources) to reach her destination without
taking into account the overall system performance.  In these systems,
the cost for each agent depends on the resources chosen as well as the
number of agents choosing the same resources. This results in a Nash
equilibrium, i.e., a point where no agent can benefit by altering its
actions, assuming that the actions of all the other agents remain fixed
\cite{Youn2008}. However, it is known that the \emph{user optimal}
policy leading to a Nash equilibrium is generally inefficient and
results in a suboptimal behavior compared to the \emph{socially optimal}
policy that could be attained through a centrally controlled system
\cite{Youn2008}. In order to quantify this inefficiency due to selfish
driving, we define the \emph{Price of Anarchy (PoA)} as the ratio of the
total travel latency cost under the \emph{user optimal (user-centric)}
routing policy vs. the \emph{socially optimal (system-centric)} one. The
PoA is, therefore, a measure of the efficiency achieved by any
transportation network as it currently operates.

The first issue addressed in the paper is how to measure the PoA from
data. The user flow equilibrium in a transportation network is known as
a Wardrop equilibrium \cite{wardrop1952some} (an instantiation of the
generic Nash equilibrium). It is the solution of the \emph{Traffic
  Assignment Problem (TAP)} \cite{patriksson1994traffic}, which we call
the \emph{user-centric forward problem}. To solve this TAP, we need to
know a priori: (\rmnum{1}) the specific travel latency cost functions
involved \cite{branston1976link}, and (\rmnum{2}) the traffic demand
expressed through an \emph{Origin-Destination (OD)} demand matrix
\cite{patriksson1994traffic}. Starting from the equilibrium link flows
(assuming they can be inferred or directly observed), we first estimate
an initial OD demand matrix. We note that the OD demand estimation
problem has been widely-studied; see, e.g.,
\cite{abrahamsson1998estimation,bera2011estimation}, and the references
therein. Then, based on inverse optimization techniques recently
developed in \cite{bertsimas2014data}, we propose a novel
\emph{user-centric inverse problem} formulation. Specifically, given
observed link flow data (Wardrop equilibrium), we estimate the
associated travel latency cost functions. In other words, we seek cost
functions which, when applied to the TAP, would yield the link flows
that are actually observed. Once this is accomplished, based on a
bi-level optimization problem formulation considered in
\cite{spiess1990gradient,lundgren2008heuristic}, we develop an
algorithmic procedure for iteratively adjusting the values of the OD
demands so that the observed link flows are as close as possible to the
solution of the \emph{user-centric forward problem} (i.e., TAP). The OD
demand and the user travel latency cost functions, completely
parametrize a predictive model of the transportation network. We use
this model to calculate the total travel latency cost under the user
optimal routing policy, thus obtaining the numerator of the PoA ratio.

Next, using the same predictive model, we formulate a
\emph{system-centric forward problem}
\cite{patriksson1994traffic,Sepid2016}, a Non-Linear Program (NLP), in
which all agents (drivers) cooperate to optimize the overall system
performance. Its solution enables us to calculate the total travel
latency cost under the socially optimal routing policy, i.e., the
denominator of the PoA ratio. Thus, the combination of the inverse and
forward optimization problems results in measuring the PoA for a given
transportation network whose equilibrium link flows are observed based
on collected traffic data.

Having an accurate predictive model allows us to go beyond estimation
(of the PoA) and consider specific control actions that could reduce the
PoA. To that end, we analyze the sensitivity of the optimal objective
function value of an optimization problem formulation for the TAP with
respect to key parameters, such as road capacities and free-flow travel
times. The results can help prioritize road segments for interventions
that can mitigate congestion. We derive sensitivity analysis formulae
and propose their finite difference approximations.

As an illustration of our data-driven approach outlined above, we use
actual traffic data from the Eastern Massachusetts (EMA) transportation
network, in the form of spatial average speeds and road segment flow
capacities. These data were provided to us by the Boston Region
Metropolitan Planning Organization (MPO) and include average speeds over
13,000 road segments at every minute of the year 2012. By using a
traffic flow model, we first infer equilibrium flows on each road
segment and then apply our approach to evaluate the PoA for two highway
subnetworks of the EMA network. In addition, we derive sensitivity
analysis results and conduct a meta-analysis comparing the user-centric
and socially optimal routing policies.

As a final step, we propose strategies for reducing the PoA. First, by
taking advantage of the rapid emergence of Connected Automated Vehicles
(CAVs) \cite{alonso2011autonomous,lee2013sustainability,zhang2016optimal,rios2017survey}, it has become feasible to automate routing decisions, thus
solving a \emph{system-centric forward problem} in which all CAVs
(bypassing driver decisions) cooperate to optimize the overall system
performance. Second, we propose a modification to existing GPS
navigation algorithms recommending to all drivers socially optimal
routes. Finally, our sensitivity analysis results provide the means to
prioritize road segments for specific interventions that can mitigate
congestion.

The rest of the paper is organized as follows. We review the related
literature in Sec.~\ref{sec:lit}. In Sec.~\ref{cdc16-sec:mod}, we
introduce models and methods we use. In Sec.~\ref{sec: dataset}, we
describe the datasets and explain the data processing procedures for a
case study of the EMA network. Numerical results for the case study are
shown in Sec.~\ref{Sec:Rsults}. In Sec.~\ref{sec:strat}, we propose
possible strategies to reduce the PoA. We provide concluding remarks and
point out some directions for future research in Sec. ~\ref{sec:conc}.

\textbf{Notation:} All vectors are column vectors. For economy of space,
we write $\bx = \left(x_1, \ldots, x_{\text{dim}(\bx)}\right)$ to denote
the column vector $\bx$, where $\text{dim}(\bx)$ is the dimension of
$\bx$. We use $\bzero$ and $\mb{1}$ for the vectors with all entries
equal to zero and one, respectively.  We denote by $\mathbb{R}_+$ the
set of all nonnegative real numbers. $\bM \geq \textbf{0}$ (resp., $\bx
\geq \textbf{0}$) indicates that all entries of a matrix $\bM$ (resp.,
vector $\bx$) are nonnegative. We use ``prime'' to denote the transpose
of a matrix or vector. Unless otherwise specified, $\|\cdot\|$ denotes
the $\ell_2$ norm. We let $\left| \mathcal{D} \right|$ denote the
cardinality of a set $\mathcal{D}$, and $\left[\kern-0.15em\left[ \scrD
  \right]\kern-0.15em\right]$ the set $\left\{ {1, \ldots ,\left| \scrD
    \right|} \right\}$.

\section{Related Work}  \label{sec:lit}

The classical static Traffic Assignment Problem (TAP)
\cite{patriksson1994traffic}, i.e., the \emph{user-centric forward
  problem} in our terminology, has been widely studied; see, e.g.,
\cite{dafermos1969traffic,leblanc1975efficient} for the
\emph{single-class} (i.e., all vehicles are modeled as belonging to the
same \emph{class}) transportation networks and
\cite{dafermos1972traffic,nagurney2000multiclass,ryu2015solving} for the
multi-class (i.e., different types of vehicles, such as cars or trucks,
are modeled as belonging to different \emph{classes}) transportation
networks. The static TAP has also been generalized to the case that has a dynamic network equilibrium modeling capability; see, e.g., \cite{friesz1989dynamic,janson1991dynamic}, among others.

Based on road traffic counts within selected time intervals (i.e., road
traffic flows), the problem of estimating the Origin-Destination (OD)
demand matrix of a given transportation network has been considered in
\cite{nguyen1984estimating,abrahamsson1998estimation,bera2011estimation},
and references therein. In particular, \cite{hazelton2000estimation}
proposed a Generalized Least Squares (GLS) method to estimate the OD
demand matrices of uncongested networks, and
\cite{spiess1990gradient,yang2001simultaneous,lundgren2008heuristic}
considered networks that could include congested roads.

Sensitivity analyses of traffic equilibria were conducted in
\cite{friesz1989dynamic,yang1997sensitivity,patriksson2004sensitivity},
among others, by evaluating the directions of change that occur in the
link flows with respect to the change of travel costs as parameters in
the cost and demand functions.

Preliminary PoA evaluation results of this paper have been presented in
two conferences, \cite{CDC16} and \cite{poa-ifac-2017}, {where results of a case study for a much smaller subnetwork of EMA were reported and no PoA reduction strategies were proposed}. A similar topic was
also discussed in \cite{monnot2017bad} and the references therein; in
particular, based on real traffic data from the transportation network
of Singapore, \cite{monnot2017bad} used a different framework from ours
to quantify the PoA.

\section{Models and Methods} \label{cdc16-sec:mod}

\subsection{Model for a single-class transportation network} \label{model}

We begin by reviewing the model of \cite{CDC16}.  Denote a road network
by $\left( {\mathcal{V}, \scrA, \mathcal{W}} \right)$, where $\left(
  {\mathcal{V}, \scrA} \right)$ forms a directed graph with
$\mathcal{V}$ being the set of nodes and $\scrA$ the set of links, and
$\mathcal{W} = \left\{ {{\bw_i}:{\bw_i} = \left( {{w_{si}},{w_{ti}}}
    \right), \,i \in [\kern-0.15em[ \scrW ]\kern-0.15em]} \right\}$
indicates the set of all OD pairs. {Note that only nodes of the road
  network can be origin/destination of flows; we make this standard
  modeling assumption to accommodate our graph-based view of the
  transportation system.} Assume the graph $\left( 
  {\mathcal{V}, \scrA} \right)$ is strongly connected and let $\bN \in
{\left\{ {0,1, - 1} \right\}^{\left| \mathcal{V} \right| \times \left|
      \scrA \right|}}$ be its node-link incidence matrix. Denote by
$\textbf{e}_{a}$ the vector with an entry being 1 corresponding to link
$a$ and all the other entries being 0.  For any OD pair $\bw = \left(
  {{w_s},{w_t}} \right)$, denote by ${d^{\bw}} \ge 0$ the amount of the
flow demand from $w_s$ to $w_t$. Let ${\bd^{\bw}} \in
{\mathbb{R}^{\left| \mathcal{V} \right|}}$ be the vector which is all
zeros, except for two entries $-{d^{\bw}}$ and ${d^{\bw}}$ corresponding
to nodes $w_s$ and $w_t$ respectively.

Denote by $\mathcal{R}_{i}$ the set of simple routes (a route without
cycles is called a ``\textit{simple route}'') for OD pair $i$.  For each
$a \in \scrA$, $i \in [\kern-0.15em[ \scrW ]\kern-0.15em]$, $r \in
\mathcal{R}_{i}$, define the link-route incidence by
\[
{\delta ^{i}_{ra}} = \left\{ \begin{gathered}
1,{\text{  if route }}r \in {\mathcal{R}_{i}}{\text{ uses link }}a, \hfill \\
0,{\text{  otherwise.}} \hfill \\
\end{gathered}  \right.
\]

Let $x_a$ denote the flow on link $a \in \scrA$ and $\bx = (x_a; \; a
\in \scrA)$ the flow vector.  Denote by ${t_a(\bx)}:\mathbb{R}_ +
^{\left| \scrA \right|} \to {\mathbb{R}_ + }$ the \textit{travel latency
  cost} (i.e., \textit{travel time}) function for link $a \in \scrA$. If
for all $a \in \scrA$, ${t_a(\bx)}$ only depends on $x_a$, we say the
cost function $\bt\left( \bx \right) = \left( {{t_a}\left( {{x_a}}
    \right); \; a \in \scrA} \right)$ is \textit{separable}
\cite{patriksson1994traffic}.  Throughout the paper, we assume that the
travel latency cost functions are separable and take the following form
\cite{bertsimas2014data,branston1976link}:
\begin{align}
{t_a}\left( {{x_a}} \right) = {t^{0}_a}f\left( {\frac{{{x_a}}}{{{m_a}}}} \right), \label{cdc16-costf}
\end{align}
where $t^{0}_a$ is the \textit{free-flow travel time} of $a \in
\scrA$, $f(0)=1$, $f(\cdot)$ is strictly increasing and continuously
differentiable on $\mathbb{R}_+$, and $m_a$ is the \textit{flow capacity} of
$a \in \scrA$. {Note that the flow capacity is not a ``hard'' constraint; $x_a$ could exceed $m_a$ for various $a$ at the cost of increased travel time.}

Define the set of feasible flow vectors $\scrF$ as \cite{bertsimas2014data}:
\[
\begin{array}{l}
\scrF \defeq \Big\{ {\bx:\exists {\bx^{\bw}} \in \mathbb{R}_ +
	^{\left| \scrA \right|} ~\text{s.t.}~\bx =
	\sum\limits_{\bw \in \scrW} {{\bx^{\bw}}},} \Big.\\
~~~~~~~~~\Big. \bN{\bx^{\bw}} = {\bd^{\bw}},\,\forall \bw \in \scrW \Big\},
\end{array}
\]
where $\bx^\bw$ indicates the flow vector attributed to OD pair $\bw$.
In order to formulate appropriate forward and inverse optimization
problems arising in transportation networks, we next state the
definition of \textit{Wardrop equilibrium}.

\begin{defi}[\cite{patriksson1994traffic}] \label{cdc16-def1} \emph{A
    feasible flow $\bx^* \in \scrF$ is a {\em{Wardrop equilibrium}} if
    for every OD pair $\bw = \left( {{w_s},{w_t}} \right) \in
    \mathcal{W}$, and any route connecting $(w_s,w_t)$ with positive
    flow in $\bx^*$, the cost of traveling along that route is no
    greater than the cost of traveling along any other route that
    connects $(w_s,w_t)$. Here, the cost of traveling along a route is
    the sum of the costs of each of its constituent links.}
\end{defi}

\subsection{The user-centric forward problem}  \label{sec:forward}

As in \cite{CDC16}, here we refer to the classical static Traffic
Assignment Problem (TAP) as the \textit{user-centric forward problem},
whose goal is to find the Wardrop equilibrium for a given single-class
transportation network with a given travel latency cost function and a
given OD demand matrix.  It is a well-known fact that, for network
$\left( {\mathcal{V}, \scrA, \mathcal{W}} \right)$, the TAP can be
formulated as the following optimization problem
\cite{dafermos1969traffic,patriksson1994traffic}:
\begin{align}
\text{(userOpt)} \quad \mathop {\min }\limits_{{\bx \in \scrF}}
\sum\limits_{a \in \scrA} {\int_0^{{x_a}} {t_a(s)ds}
}. \label{cdc16-tap} 
\end{align}
As an alternative, we also formulate the TAP as a Variational Inequality
(VI) problem:

\begin{defi} [\cite{bertsimas2014data}]
	\label{cdc16-def2}  \em{The VI problem, denoted as
          $\text{VI}\left( {\bt,\scrF} \right)$, is to find an ${\bx^ *
          } \in \scrF$ s.t. 
	\begin{align}
	\bt{\left( {{\bx^ * }} \right)'}\left( {\bx - {\bx^ * }} \right)
        \geq 0, \quad \forall \bx \in \scrF. \label{cdc16-VI} 
	\end{align}}
\end{defi}

To proceed, let us first recall the definition of the strong
monotonicity for a cost function: $\bt(\cdot)$ is \textit{strongly
  monotone} \cite{patriksson1994traffic} on $\scrF$ if there exists a
constant $\eta > 0$ such that
\begin{align}\left[ {\bt\left( \bx \right) - \bt\left( \by
              \right)} \right]'\left( {\bx - \by} \right) \ge \eta
          {\left\| {\bx - \by} \right\|^2},\quad \forall \bx,\by \in
          \scrF.   \label{str-mon} 
\end{align}
It is known that if $\bt(\cdot)$ is continuously differentiable on
$\scrF$, then \eqref{str-mon} is equivalent to the positive definiteness
of the Jacobian of $\bt(\cdot)$
\cite[p. 180]{patriksson1994traffic}. Note that a strictly increasing
$f(\cdot)$ in \eqref{cdc16-costf} would not necessarily ensure the
strong monotonicity of $\bt(\cdot)$; e.g., $f(x) \defeq x^3$ and
$\bt(\bx) \defeq ({x_1^3}, {x_2^3})$ would lead to the Jacobian of
$\bt(\bx)$ as 
\[ \begin{bmatrix}
      3x_1^2 & 0 \\
      0 & 3x_2^2
	\end{bmatrix},
\]
which is not positive definite over $\mathbb{R}^2$.
We next introduce a key assumption.

\begin{ass}
	\label{cdc16-assumption1}
	\emph{$\bt(\cdot)$ is strongly monotone on $\scrF$ and
          continuously differentiable on $\mathbb{R}_ + ^{\left| \scrA
            \right|}$. $\scrF$ is nonempty and contains an interior
          point (Slater's condition \cite{boyd2004convex}).}
\end{ass}
For the existence and uniqueness of the TAP, the following result is
available:
\begin{thm}[\cite{patriksson1994traffic}] \label{cdc16-th21}
  \textit{Assumption~\ref{cdc16-assumption1}} implies that there exists a
  Wardrop equilibrium of the network $\left( {\mathcal{V}, \scrA,
      \mathcal{W}} \right)$, which is the unique solution to
  $\text{VI}(\bt,\scrF)$.
\end{thm}

\subsection{The user-centric inverse problem} \label{cdc16-sec InverseVI-uni}

To solve the \emph{user-centric forward problem}, we need to know the
travel latency cost function and the OD demand matrix. Assuming that we
know the OD demand matrix and have observed the Wardrop equilibrium link
flows, we seek to formulate the \emph{user-centric inverse problem} (the
inverse VI problem, in particular), so as to estimate the travel latency
cost function.  To provide some insight, given $|\scrK|$ samples of the
link flow vector $\bx$, one can think of them as flow observations on
$|\scrK|$ different networks/subnetworks which are nevertheless produced
by the exact same cost function. The inverse formulation seeks to
determine the cost function so that each flow observation is as close to
an equilibrium as possible. Given that the inverse problem will rely on
measured flows, we should expect measurement noise which will prevent
the flows from being an exact solution of the forward VI problem
$\VI(\mathbf{t}, \scrF)$. Therefore, we will first define the notion of
an approximate solution.

For a given $\epsilon > 0$, we define an
\textit{$\epsilon$-approximate solution} to $\VI(\mathbf{t}, \scrF)$ by
changing the right-hand side of \eqref{cdc16-VI} to $- \epsilon$:

\begin{defi}[\cite{bertsimas2014data}]
	\label{cdc16-def3}
	\emph{Given $\epsilon > 0$,
	$\hat{\bx} \in \scrF$ is called an
	\textit{$\epsilon$-approximate solution} to $\VI(\mathbf{t}, \scrF)$ if
	\begin{equation}
	\label{cdc16-eq:defApproxEquil}
	\mathbf{t}(\hat{\bx})'(\bx - \hat{\bx}) \geq - \epsilon, \quad
	\forall \bx \in \scrF.
	\end{equation}}
\end{defi}

Assume now we are given $|\scrK|$ networks
$(\scrV^{\scriptscriptstyle{(k)}},\scrA^{\scriptscriptstyle{(k)}},\scrW^{\scriptscriptstyle{(k)}}),\,
k \in [\kern-0.15em[ \scrK ]\kern-0.15em]$ (as a special case, these
could be $|\scrK|$ replicas of the same network $\left( {\scrV, \scrA,
    \scrW} \right)$), and the observed link flow data
$\{\bx^{\scriptscriptstyle{(k)}} = (x_a^{\scriptscriptstyle{(k)}}; \, a
\in \scrA^{\scriptscriptstyle{(k)}}); \, k \in [\kern-0.15em[ \scrK
]\kern-0.15em] \}$, where $k$ is the network index and
$x_a^{\scriptscriptstyle{(k)}}$ is the flow for link $a \in
\scrA^{\scriptscriptstyle{(k)}}$ correspondingly. The inverse VI problem
amounts to finding a function $\bt$ such that
$\bx^{\scriptscriptstyle{(k)}}$ is an $\epsilon_k$-approximate solution
to $\VI(\bt, \scrF^{\scriptscriptstyle{(k)}})$ for each $k$.  Denoting
$\bepsilon = (\epsilon_k; \, k \in [\kern-0.15em[ \scrK
]\kern-0.15em])$, we can formulate the inverse VI problem as
\cite{bertsimas2014data}:
\begin{align}
\min_{{\mathbf{t}}, \bepsilon} \quad & \| \bepsilon \|
\label{cdc16-inverseVI1} \\
\text{s.t.} \quad & {\mb{t}}({\bx^{\scriptscriptstyle{(k)}}})'(\bx-{\bx^{\scriptscriptstyle{(k)}}}) \geq -\epsilon_k, \quad \forall
\bx\in \scrF^{\scriptscriptstyle{(k)}}, \,k \in [\kern-0.15em[ \scrK 
]\kern-0.15em], \notag \\
& \epsilon_k > 0, \quad \forall k \in [\kern-0.15em[ \scrK 
]\kern-0.15em], \notag
\end{align}
where the optimization is over the selection of function $\mathbf{t}$
and the vector $\bepsilon$.

Aiming at recovering a cost function $\bt$ that has both good data
reconciling and generalization properties (i.e., $\bt$ should fit
``old'' data well but should not be overfitting; it must also have great
power to predict ``new'' data), to make \eqref{cdc16-inverseVI1}
solvable, we apply an estimation approach which expresses the function
$f(\cdot)$ (in \eqref{cdc16-costf}) in a Reproducing Kernel Hilbert
Space (RKHS) $\scrH$
\cite{bertsimas2014data,evgeniou2000regularization}. In particular, by
\cite[Thm. 2]{bertsimas2014data}, we reformulate the inverse VI problem
\eqref{cdc16-inverseVI1} as
	\begin{align} 
	&\text{(invVI-1)} \quad 
	\min_{{f}, \by, \bepsilon} \quad    \left\| \boldsymbol{\epsilon}  \right\| + \gamma\| {f} \|^2_{\scrH}  \label{nonPara}
	\\   
	\text{s.t.} \quad & 
	{{\textbf{e}}}_a'\bN_k'{\by^\bw} \leq t_a^0{{f} {\left( {\frac{{{x_a}}}{{{m_a}}}} \right)}}, \; \label{dual-feasi} 
	\\ &~~~~~~~~~~ \forall \bw \in {\mathcal{W}^{\scriptscriptstyle{(k)}}}, ~a \in {\scrA^{\scriptscriptstyle{(k)}}},~k \in [\kern-0.15em[ \scrK 
	]\kern-0.15em], \notag \\ 
	&\sum\limits_{a \in {\scrA^{\scriptscriptstyle{(k)}}}}{t_a^0{x_a}{ {{f}\left( {\frac{{{x_a}}}{{{m_a}}}} \right)} }} - \sum\limits_{\bw \in {\mathcal{W}^{\scriptscriptstyle{(k)}}}} {\left( {{\bd^\bw}} \right)'{\by^\bw}}  \leq {\epsilon _k},\label{primal-dual-gap}
	\\ &~~~~~~~~~~ \forall k \in [\kern-0.15em[ \scrK 
	]\kern-0.15em], \notag \\
	&{ {{f}\left( {\frac{{{x_a}}}{{{m_a}}}} \right)} } \leq { {{f}\left( {\frac{{{x_{\tilde a}}}}{{{m_{\tilde a}}}}} \right)} }, \label{increase}
	\\ &~~~~~~~~~~ \forall a,~\tilde a \in {{ \bigcup\nolimits_{k =
				1}^{|\scrK|} {{\scrA^{\scriptscriptstyle{(k)}}}}}}~\text{s.t. }\frac{{{x_a}}}{{{m_a}}} \leq \frac{{{x_{\tilde a}}}}{{{m_{\tilde a}}}}, \notag \\
		& \bepsilon \geq \textbf{0},  \quad 
		{f} \in \scrH,  \notag
		\\  
	&{f}(0) = 1, \label{normaliz}
\end{align}
which is a counterpart of \cite[(22)]{bertsimas2014data}. Note that $\by
= \left( {{\by^{\bw}};\,\bw \in {\scrW^{\scriptscriptstyle{(k)}}},k \in
    [\kern-0.15em[ \scrK ]\kern-0.15em]} \right)$ and
$\boldsymbol{\epsilon} = (\epsilon_k; \,k \in [\kern-0.15em[ \scrK
]\kern-0.15em])$ are decision vectors ($\by^{\bw}$ is a dual variable
which can be interpreted as the ``price'' of $\bd^{\bw}$, in
particular). Note also that $\gamma > 0$ is a regularization parameter
-- a smaller $\gamma$ should result in recovering a ``tighter''
$f(\cdot)$ in terms of data reconciling; a larger $\gamma$, on the other
hand, would lead to a ``better'' $f(\cdot)$ in terms of generalization
properties. Moreover, $\| {f} \|^2_{\scrH}$ denotes the squared norm of
$f(\cdot)$ in $\scrH$, \eqref{dual-feasi} is for dual feasibility,
\eqref{primal-dual-gap} is the suboptimality (primal-dual gap)
constraint, \eqref{increase} enforces $f(\cdot)$ to be non-decreasing,
and \eqref{normaliz} is a normalization constraint.
	
It can be seen that the above formulation is still too abstract for us
to solve, because it is an optimization over functions. To make it
tractable, in the following we will specify $\scrH$ by selecting its
\textit{reproducing kernel} \cite{evgeniou2000regularization} to be a
polynomial ${\phi}(x, y) = ( c + xy )^n$ for some choice of $c \geq 0$
and $n \in \mathbb{N}$ (for the specifications of $c$ and $n$, see
\cite{InverseVIsTraffic}). Then, writing
\[\phi\left( {x,y} \right) = {\left( {c + xy} \right)^n} =
\sum\limits_{i = 0}^n {{n \choose i}{c^{n - i}}{x^i}{y^i}}, \] by
\cite[(3.2), (3.3), and (3.6)]{evgeniou2000regularization}, we
instantiate invVI-1 as
\begin{align}
&\text{(invVI-2)} \quad
\mathop {\min }\limits_{\bbeta ,\by,\boldsymbol{\epsilon} } \quad  \left\| \boldsymbol{\epsilon}  \right\| + \gamma \sum\limits_{i = 0}^n {\frac{{\beta _i^2}}{{{n \choose i}{c^{n - i}}}}} \notag \\
\text{s.t.}\quad &{\textbf{e}}_a'\bN_k'{\by^\bw} \leq t_a^0{\sum\limits_{i = 0}^n {{\beta _i}\left( {\frac{{{x_a}}}{{{m_a}}}} \right)} ^i}, \; \notag \\ 
&~~~~~~~~~~ \forall \bw \in {\mathcal{W}^{\scriptscriptstyle{(k)}}}, ~a \in {\scrA^{\scriptscriptstyle{(k)}}},~k \in [\kern-0.15em[ \scrK 
]\kern-0.15em], \notag \\ 
&\sum\limits_{a \in {\scrA_k}}{t_a^0{x_a}{\sum\limits_{i = 0}^n {{\beta _i}\left( {\frac{{{x_a}}}{{{m_a}}}} \right)} ^i}} - \sum\limits_{\bw \in {\mathcal{W}_k}} {\left( {{\bd^\bw}} \right)'{\by^\bw}}  \leq {\epsilon _k},\;\notag \\ &~~~~~~~~~~ \forall k \in [\kern-0.15em[ \scrK 
]\kern-0.15em], \notag \\
&{\sum\limits_{i = 0}^n {{\beta _i}\left( {\frac{{{x_a}}}{{{m_a}}}} \right)} ^i} \leq {\sum\limits_{i = 0}^n {{\beta _i}\left( {\frac{{{x_{\tilde a}}}}{{{m_{\tilde a}}}}} \right)} ^i},\;\notag \\ &~~~~~~~~~~ \forall a,~\tilde a \in {{ \bigcup\nolimits_{k =
			1}^{|\scrK|} {{\scrA^{\scriptscriptstyle{(k)}}}}}}~\text{s.t. }\frac{{{x_a}}}{{{m_a}}} \leq \frac{{{x_{\tilde a}}}}{{{m_{\tilde a}}}}, \notag \\
&\bepsilon \ge \textbf{0}, \quad {\beta _0} = 1, \notag 
\end{align}
where the function $f(\cdot)$ in invVI-1 is parameterized by 
$\bbeta  = \left( {{\beta _i};\;i = 0, 1, \ldots ,n} \right)$.
Assuming an optimal ${\bbeta ^*} = \left( {\beta _i^*;\,i = 0,1, \ldots ,n} \right)$ is obtained by solving invVI-2, then our estimator for $f(\cdot)$ is
\begin{align}
\hat f\left( x \right) = \sum\limits_{i = 0}^n {\beta _i^*{x^i}}  = 1 + \sum\limits_{i = 1}^n {\beta _i^*{x^i}}.  \label{cdc16-costEstimator}
\end{align}

\subsection{OD demand estimation}

Given a network $\left( {\mathcal{V}, \scrA, \mathcal{W}} \right)$, to
estimate an initial OD demand matrix, we borrow the Generalized Least
Squares (GLS) method proposed in \cite{hazelton2000estimation}, which
assumes that the transportation network $\left( {\mathcal{V}, \scrA,
    \mathcal{W}} \right)$ is uncongested (in other words, for each OD
pair the route choice probabilities are independent of traffic flow),
and that the OD trips (traffic counts) are Poisson distributed. Note
that such assumptions may be strong and we will relax them when
finalizing our OD demand estimator by performing an adjustment
procedure.

Denote by $\{ {{\bx^{\scriptscriptstyle{(k)}}}}; \, k \in [\kern-0.15em[
\scrK ]\kern-0.15em]\}$ $|\scrK|$ observations of the flow vector. Let
$\bar \bx = (1/|\scrK|)\sum\nolimits_{k = 1}^{|\scrK|}
{{\bx^{\scriptscriptstyle{(k)}}}} $ be the sample mean vector and
\[ \bS =
\frac{1}{(|\scrK| - 1)}\sum_{k = 1}^{|\scrK|} {\big(
  {{\bx^{\scriptscriptstyle{(k)}}} - \bar {\bx}} \big){{\big(
      {{\bx^{\scriptscriptstyle{(k)}}} - \bar {\bx}} \big)}' }}
\]
 the
sample covariance matrix. Let $\bP = \left[ {{p_{ir}}} \right]$ denote
the route choice probability matrix, where $p_{ir}$ is the probability
that a traveler associated with OD pair $i$ uses route $r$. Vectorize
the OD demand matrix as $\bg = (g_i; \,i \in [\kern-0.15em[ \scrW
]\kern-0.15em]) $. After finding feasible routes for each OD pair, thus
obtaining the link-route incidence matrix $\bA$, the GLS method amounts
to solving the following optimization problem:
\begin{align}
\text{(P0)} \quad \mathop {\min }\limits_{\bP \geq \textbf{0}, \,\bg
  \geq \textbf{0} } \quad & \sum\limits_{k = 1}^{|\scrK|} {\left(
    {{\bx^{\left( k \right)}} - \bA\bP'\bg } \right)'{\bS^{ - 1}}\left(
    {{\bx^{\left( k \right)}} - \bA\bP'\bg } \right)}  \notag \\ 
\text{s.t.} \quad 
& p_{ir} = 0 \quad \forall (i,r) \in \{(i,r): r \notin \scrR_i\}, \notag \\
&\bP\textbf{1} = \textbf{1}, \notag 
\end{align}
which minimizes a weighted sum of the squared errors in the flow
observations. Directly solving (P0) is cumbersome due to the complicated
form of the objective function, and we in turn decouple (P0) into two
subproblems. To that end, we perform a variable substitution by setting
$\bxi = \bP'\bg$ and we let $h(\bP, \bg)$ be an arbitrarily selected smooth
scalar-valued function. Then, we solve sequentially the following two
problems \cite{CDC16}:
\begin{align}
\text{(P1)} \quad \mathop {\min }\limits_{\bxi  \geq \textbf{0}} \quad
\frac{|\scrK|}{2}\bxi '\bQ\bxi  - \bb'\bxi,   \label{qp2} 
\end{align} 
where $\bQ = \bA'{\bS^{ - 1}}\bA$ and $\bb =
\sum\nolimits_{k = 1}^{|\scrK|} {\bA'{\bS^{ -
      1}}{\bx^{\scriptscriptstyle{(k)}}}} $, and 
\begin{align}
\text{(P2)} \quad &\mathop {\min }\limits_{\bP \geq \textbf{0}, \,\bg
  \geq \textbf{0} } \quad h\left( {\bP,\bg } \right) \label{qp3} \\ 
\text{s.t.} \quad  
& p_{ir} = 0, \quad \forall (i,r) \in \{(i,r): r \notin \scrR_i\}, \notag \\
&\bP'\bg  = {\bxi ^0}, \notag \\
&\bP\textbf{1} = \textbf{1}, \notag 
\end{align}
where $\bxi^0$ is the optimal solution to (P1). Essentially, (P1) uses
the variable substitution to eliminate the constrains on $\bP$ and (P2)
seeks to find a feasible $\bP$ consistent with the optimal solution of
(P1) and the relationship $\bxi = \bP'\bg$. We write the feasibility
problem (P2) as an optimization problem with some ``dummy'' cost
function because this allows us to use an optimization solver; in fact,
we can simply set $h(\bP,\bg) \equiv 0$. Specifically, (P1) (resp.,
(P2)) is a typical \textit{Quadratic Program (QP)} (resp.,
\textit{Quadratically Constrained Program (QCP)}). Letting $(\bP^0,
\bg^0)$ be an optimal solution to (P2), then $\bg^0$ is our initial
estimate of the demand vector.

\begin{rmk} \label{rem:od}
\em{It is seen that each entry of $\bg^0$ can always be expressed as a sum of
  certain entries in $\bxi^0$; in other words, given $\bxi^0 \ge \textbf{0}$,
  (P2) always has a feasible solution. Thus, (P0) is actually equivalent to
  (P1) and (P2), in the sense that if $(\bP^0,\bg^0)$ is an optimal solution
  to (P0) (resp., (P2)), then it is also an optimal solution to (P2) (resp.,
  (P0)).  In addition, we note that the GLS method above would encounter
  numerical difficulties when the network size is large, because there would
  be too many decision variables. Note also that this method is valid under a
  ``no-congestion'' assumption and, to take the congestion on the link flows
  into account, we in turn consider a bi-level optimization problem in the
  following.}
\end{rmk}

Assume now the function $f(\cdot)$ in \eqref{cdc16-costf} is
available. For any given feasible $\bg$ ($\geq \textbf{0}$), let
$\bx(\bg)$ be the optimal solution to the TAP \eqref{cdc16-tap}.  In the
following, denote by $\tilde{\bx} = (\tilde{x}_a; \;a \in \scrA)$ the
observed flow vector. Assuming an initial demand vector $\bg^{0}$ is
given (the $\bg^0$ obtained by solving (P1) and (P2) is a good
candidate; we will take it as $\bg^{0}$ hereafter), we consider the
following bi-level optimization problem
\cite{spiess1990gradient,lundgren2008heuristic}:

{\begin{align}{\text{(BiLev)}}\;\;\;\mathop {\min }\limits_{\bg \geq \textbf{0}} {\text{  }}F\left( \bg \right) \defeq & {\gamma _1}\sum\limits_{i \in [\kern-0.15em[ \scrW 
		]\kern-0.15em]} {{{\left( {{g_i} - g_i^0} \right)}^2}} \notag \\
	& + {\gamma _2}\sum\limits_{a \in \scrA} {{{\left( {{x_a}\left( \bg \right) - {{\tilde x}_a}} \right)}^2}}, \label{bilevel}
\end{align}
where $\gamma_1, \gamma_2 \geq 0$ are two weight parameters. The first
term penalizes moving too far away from the initial demand, and the
second term ensures that the optimal solution to the TAP is close to the
flow observation.} Note that the BiLev formulation \eqref{bilevel} is
more general than the one considered in \cite{poa-ifac-2017}, {which includes the second term only}.  It is worth
pointing out that $F\left( \bg \right)$ has a lower bound $0$ which
guarantees the convergence of the algorithm (see
Alg. \ref{alg:demandAdjustment}) that we will apply.

\begin{rmk} \label{rem:gamma1}
\em{From now on, let us fix $\gamma_2=1$ in \eqref{bilevel}. Intuitively, the
  closer the initial $\bg^0$ to the ground truth $\bg^*$, the larger the
  $\gamma_1$ we should set; otherwise the contribution of the first term to
  the objective function will be small. In practice, however, we typically do
  not have exact information about how far $\bg^0$ is from $\bg^*$; we
  therefore have to appropriately tune $\gamma_1$. One possible criterion
  is that, fixing the parameters involved in Alg. \ref{alg:demandAdjustment},
  a ``good'' $\gamma_1$ should lead to a reduction of the objective function
  value of the BiLev as much as possible.  }
\end{rmk}

To solve the BiLev numerically, thus adjusting the demand vector
iteratively, we leverage a gradient-based algorithm
(Alg. \ref{alg:demandAdjustment}). In particular, suppose that the route
probabilities are locally constant. For OD pair $i \in [\kern-0.15em[
\scrW ]\kern-0.15em]$, we consider only the {\emph{fastest}} route
$r_i(\bg)$, where in each iteration, based on the updated link
  flows after the previous iteration, we update link \emph{travel times}
  and assign them as current link weights in the graph model introduced
  in Sec. \ref{model}.  Then, we have \cite{spiess1990gradient}
\begin{align}
  \frac{{\partial {x_a}\left( \bg \right)}}{{\partial {g _i}}} \approx
  {\delta _{{r_i(\bg)}a}} = \left\{ \begin{gathered}
      1,{\text{ if }}a \in {r_i(\bg)}, \hfill \\
      0,{\text{ otherwise}}{\text{.}} \hfill \\
\end{gathered}  \right. \label{jacobbb}
\end{align}
(Note that we have assumed the partial derivatives do exist; a
comprehensive discussion on the existence and calculation of $\partial
{x_a}(\bg)/\partial g _i$ can be found in
\cite{patriksson2004sensitivity}.)  Thus, by \eqref{jacobbb} we obtain
an approximation to the Jacobian matrix
\begin{align}
\left[ {\frac{{\partial {x_a}\left( {{\bg}} \right)}}{{\partial {g _i}}};\;a \in \scrA,\;i \in [\kern-0.15em[ \scrW 
	]\kern-0.15em]} \right]. \label{jacoMat}
\end{align}
Let us now compute the gradient of $F\left( \bg \right)$. We have 
\begin{align}
&\nabla F\left( \bg \right) = \left( {\frac{{\partial F\left( \bg \right)}}{{\partial {g_i}}};\;i \in [\kern-0.15em[ \scrW 
	]\kern-0.15em]} \right) \notag \\
&= \bigg(2{\gamma _1}\left( {{g_i} - g_i^0} \right) + 2{\gamma _2}\sum\limits_{a \in \scrA} {\left( {{x_a}\left( \bg \right) - {{\tilde x}_a}} \right)\frac{{\partial {x_a}\left( \bg \right)}}{{\partial {g_i}}}} ;\; \notag \\
&~~~~~~i \in [\kern-0.15em[ \scrW 
]\kern-0.15em] \bigg).  \label{gradi}
\end{align}

\begin{rmk} \label{rem:shortest} \em{There are three reasons why we consider only the {\emph{fastest}} routes for the purpose of calculating the
    Jacobian: (\rmnum{1}) GPS navigation is widely-used by vehicle
    drivers so that they tend to always select the {\emph{fastest}} routes
    between their OD pairs. (\rmnum{2}) There are very efficient
    algorithms for finding the {\emph{fastest}} route for each OD
    pair. (\rmnum{3}) If considering more than one route for an OD pair,
    then the route flows cannot be uniquely determined by solving the
    TAP \eqref{cdc16-tap}, thus leading to unstable route-choice
    probabilities, which would undermine the accuracy of the
    approximation to the Jacobian matrix in \eqref{jacoMat}.}
\end{rmk}
We summarize the procedures for adjusting the OD demand matrices as
Alg. \ref{alg:demandAdjustment}, whose convergence will be proven in the
following proposition. 
\begin{algorithm}
	\caption{Adjusting OD demand matrices}
	\label{alg:demandAdjustment}
	\begin{algorithmic}[1]
		\Require  the road network $\left( {\mathcal{V}, \scrA, \mathcal{W}} \right)$; the function $f(\cdot)$ in \eqref{cdc16-costf}; the observed flow vector from given data $\tilde{\bx} = (\tilde{x}_a; \;a \in \scrA)$; the initial demand vector ${\bg^{0}} = (g^{0}_i; \;i \in [\kern-0.15em[ \scrW 
		]\kern-0.15em])$; two positive integer parameters $\rho$, $T$; two real parameters
		$\varepsilon_1 \geq 0$, $\varepsilon_2 > 0$. 
		\State \textbf{Step 1:} Initialization. Take the demand vector $\bg^{0}$ as the input, solve the TAP \eqref{cdc16-tap} (using the Method of Successive Averages (MSA) \cite{noriega2007algorithmic}) to obtain $\bx^0$. Set $l=0$. If $F\left(\bg^{0}\right) = 0$, stop; otherwise, go onto Step 2.
		\State \textbf{Step 2:} Computation of a descent direction. Calculate ${\bh^l} =  - \nabla F\left( {{\bg ^l}} \right)$ by \eqref{gradi}.
		\State \textbf{Step 3:} Calculation of a search direction. For $i \in [\kern-0.15em[ \scrW 
		]\kern-0.15em]$ set 
		$$\bar h_i^l = \left\{ \begin{gathered}
		h_i^l,{\text{  if }}\left( {g _i^l > \varepsilon_1 } \right){\text{ or }}\left( {g _i^l \leq \varepsilon_1 {\text{ and }}h_i^l > 0} \right){\text{, }} \hfill \\
		0,{\text{  otherwise.}} \hfill \\ 
		\end{gathered}  \right.$$
		
		\State \textbf{Step 4:} Armijo-type line search. 
		\begin{itemize}
\item[] \textbf{4.1:} Calculate the maximum possible step-size $\theta _{\max
}^l = \min \left\{ { - {{g _i^l}}\big/{{\bar h_i^l}}; \; \bar h_i^l < 0, i \in
  [\kern-0.15em[ \scrW ]\kern-0.15em]} \right\}$.
\item[] \textbf{4.2:} Determine ${\theta ^l} = \mathop {\arg \min
}\limits_{\theta \in \scrS} F\left( {{\bg ^l} + \theta {{\bar {\bh}}^l}}
  \right)$, where $\scrS \defeq \left\{{\theta _{\max }^l,{{\theta _{\max }^l}
      \mathord{\left/ {\vphantom {{\theta _{\max }^l} \rho }} \right.
        \kern-\nulldelimiterspace} \rho },{{\theta _{\max }^l} \mathord{\left/
        {\vphantom {{\theta _{\max }^l} {{\rho ^2}}}} \right.
        \kern-\nulldelimiterspace} {{\rho ^2}}}, \ldots ,{{\theta _{\max }^l}
      \mathord{\left/ {\vphantom {{\theta _{\max }^l} {{\rho ^{T - 1}}}}}
        \right.  \kern-\nulldelimiterspace} {{\rho ^{T}}}}}, 0 \right\}$.
		\end{itemize}
		
\State \textbf{Step 5:} Update and termination.
\begin{itemize}
\item[] \textbf{5.1:} Set ${\bg ^{l + 1}} = {\bg ^l} + {\theta ^l}{{\bar
    \bh}^l}$. Using $\bg^{{l+1}}$ as the input, solve the TAP
  \eqref{cdc16-tap} to obtain $\bx^{{l+1}}$.
\item[] \textbf{5.2:} If $\frac{{F\left( {{\bg^l}} \right) - F\left( {{\bg^{l
          + 1}}} \right)}}{{F\left( {{\bg^0}} \right)}} < {\varepsilon _2}$,
  stop the iteration; otherwise, go onto Step 5.3.  
\item[] \textbf{5.3:} Set $l=l+1$ and return to Step 2.	
		\end{itemize}
	\end{algorithmic}
\end{algorithm}

\begin{prop} \label{prop1} Alg. \ref{alg:demandAdjustment} converges.
\end{prop}
\pf If the initial demand vector $\bg^{0}$ satisfies $F\left(\bg^{0}\right) = 0$,
then, by Step 1, the algorithm stops (trivial case). Otherwise, we have
$F\left(\bg^{0}\right) > 0$, and it is seen from \eqref{bilevel} that the
objective function $F\left(\bg\right)$ has a lower bound $0$. In addition, by
the line search and the update steps (Steps 4.2 and 5.1, in particular), we
obtain
\begin{align}
F\left( {{\bg^{l + 1}}} \right) = & F\left( {{\bg^l} + {\theta ^l}{{\bar
      \bh}^l}} \right) \notag\\ = & \mathop {\min }\limits_{\theta \in \scrS}
F\left( {{\bg^l} + \theta {{\bar \bh}^l}} \right) \leq F\left( {{\bg^l}}
\right), \forall l, \notag
\end{align}
where the last inequality holds due to $0 \in \scrS$, indicating that the
nonnegative objective function in \eqref{bilevel} is non-increasing as the
number of iterations increases. Thus, by the well-known monotone convergence
theorem, the convergence of the algorithm can be guaranteed.
\qed

\begin{rmk}
  \em{Alg. \ref{alg:demandAdjustment} is a variant of the algorithms
    proposed in \cite{spiess1990gradient} and
    \cite{lundgren2008heuristic}. We use a different method to calculate
    the \textit{step-sizes} (resp., \textit{Jacobian matrix}) than that
    in \cite{spiess1990gradient} (resp.,
    \cite{lundgren2008heuristic}). The optimization problem BiLev is not
    convex because of the potential nonlinearity in $\bx(\bg)$. Thus,
    one would not necessarily expect Alg. \ref{alg:demandAdjustment}'s
    convergence to a global minimum. In addition, due to inaccuracies in
    the gradient calculation, one would not expect
    Alg. \ref{alg:demandAdjustment}'s convergence to a local minimum
    either. A discussion on the performance of similar heuristics can be
    found in \cite{lundgren2008heuristic}. It is worth noting that, in
    \cite{lundgren2008heuristic}, the proposed ``descent'' algorithm
    could possibly not ``descend'' in some iterations due to
    computational inaccuracy of the gradient. We will demonstrate our
    findings for the performance of Alg.~\ref{alg:demandAdjustment} by
    numerical experiments in Sec. \ref{sec:od-adj}. We also note that,
    in terms of decreasing the objective function value of the BiLev,
    the performance of Alg. \ref{alg:demandAdjustment} definitely
    depends heavily on the initial demand vector $\bg^0$.  }
\end{rmk}

\subsection{Price of Anarchy} \label{sec:poa}

As discussed in Sec. \ref{sec:intro}, one of our goals is to measure
inefficiency in the network due to the non-cooperative behavior
of drivers. Thus, we compare the network performance under a
user-centric routing policy vs. a system-centric one. As a metric for
this comparison, we conceptually define the PoA as the ratio between the
total travel latency cost, i.e., the total travel time over all drivers,
obtained under Wardrop flows (user-centric routing policy) and that
obtained under socially optimal flows (system-centric routing policy).

Given road network $\left( {\mathcal{V}, \scrA, \mathcal{W}} \right)$,
as in \cite{CDC16}, we calculate its \textit{total travel latency cost} as
\begin{align}
L(\bx) =\sum_{a \in \scrA}x_a t_a(x_a). \label{TotalLatency2}
\end{align}
The socially optimal flow vector, denoted by $\bx^{\text{social}} =
(x_a^{\text{social}};\,a \in \scrA)$, is the solution to the following
\emph{system-centric forward problem}, which is a Non-Linear Program
(NLP) \cite{patriksson1994traffic,Sepid2016}:
\begin{align}
\text{(socialOpt)} \quad \mathop {\min }\limits_{{\bx \in \scrF}} {\text{ }} 
\sum_{a \in \scrA}x_a t_a(x_a).  \label{obj_soc}
\end{align}
We therefore explicitly define the \textit{Price of Anarchy} as
\begin{equation}
\text{PoA} \defeq \dfrac{L(\bx^{\text{user}})}{L(\bx^{\text{social}})} = \frac{\sum_{a \in
		\scrA}x_a^{\text{user}}t_a(x_a^{\text{user}})}{\sum_{a \in
		\scrA}x_a^{\text{social}}t_a(x_a^{\text{social}})} \ge 1, \label{poa-def}
\end{equation}
where $\bx^{{\text{user}}} = (x_a^{{\text{user}}}; \, a \in \scrA)$ is
the Wardrop equilibrium flow vector assumed to be directly observable or
indirectly inferrable. By the definition of $\bx^{\text{social}}$, we
always have $\text{PoA} \ge 1$; the larger the PoA, the larger the
inefficiency induced by selfish drivers. Thus, PoA quantifies the
inefficiency that a societal group has to deal with due to
non-cooperative behavior of its members.

We note that the objective function in \eqref{obj_soc} is different from
its counterpart in \eqref{cdc16-tap}; for a detailed explanation, see
\cite{dafermos1969traffic}. However, the two forward problem
formulations have a very tight connection. Let us take a close look at
the following equalities \cite{patriksson1994traffic}:
\begin{align} {\overline t _a}\left( {{x_a}} \right) \defeq
  \frac{d}{{d{x_a}}}\left( {{x_a}{t_a}\left( {{x_a}} \right)} \right) =
  {t_a}\left( {{x_a}} \right) + {x_a}{\dot t_a}\left( {{x_a}} \right),\
  \forall a \in \scrA.  \label{user-social}
\end{align}
By \eqref{user-social} we see that the socialOpt in \eqref{obj_soc} is
equivalent to
\begin{align} 
\text{(userOpt)} \quad \mathop {\min }\limits_{{\bx \in \scrF}} \sum\limits_{a \in \scrA} {\int_0^{{x_a}} {\overline t_a(s)ds} }. \notag
\end{align}
The remarkable implication of the above is that in order to find the
socially optimal flows $x_a^{\text{social}}$, $a \in \scrA$, instead of
directly solving \eqref{obj_soc}, it suffices to solve \eqref{cdc16-tap}
with $t_a(\cdot)$ replaced by $\overline t_a(\cdot)$. As noted in
\cite{patriksson1994traffic}, the difference between the social cost and
the user cost is ${x_a}{\dot t_a}\left( {{x_a}} \right)$, which can be
interpreted as the cost a user (driver) imposes on the other users.

Let $\overline \bt(\bx) \defeq (\overline t_a(x_a); \, a \in \scrA)$.
To ensure the existence and uniqueness of the solution to
\eqref{obj_soc}, we need the following assumption:
\begin{ass}
	\label{cdc16-assumption2}
	\emph{$\overline \bt(\cdot)$ is strongly monotone on $\scrF$ and
	continuously differentiable on $\mathbb{R}_ + ^{\left| \scrA
		\right|}$. $\scrF$ satisfies Slater's condition
              \cite{boyd2004convex}.} 
\end{ass}
We note that if Assump. \ref{cdc16-assumption1} holds and, for all $a
\in \scrA$, $t_a(x_a)$ is convex and twice continuously differentiable
on $\mathbb{R}_ +$ (e.g., $t_a(x_a) = 2x_a^2 + x_a + 1$), then
Assump. \ref{cdc16-assumption2} holds as well.

\subsection{Sensitivity analysis}

To prioritize road segments for potential congestion reducing
interventions by the local transportation authorities, we investigate
the sensitivities of the optimal objective function value of
\eqref{cdc16-tap} with respect to key parameters, specifically,
\emph{free-flow travel time} and \emph{flow capacity}. In particular, we
first derive two rigorous formulae, and then propose their \emph{finite
  difference approximations} as an alternative.

Write ${\bt^{0}} \defeq \left(
{{t^{0}_a; \, a \in \scrA}} \right)$, $\bm
\defeq \left( {{m_a; \, a \in \scrA}} \right)$, and
\begin{align}
V\left( {{\bt^{0}},{\bm}} \right) \defeq \mathop {\min }\limits_{{\bx \in \scrF}} \sum\limits_{a \in \scrA} {\int_0^{{x_a}} {{t^{0}_a}f\left( {\frac{s}{{{m_a}}}} \right)ds} }.  \label{vtm}
\end{align}
Differentiating \eqref{vtm}, for each $a \in \scrA$ we obtain
\begin{align}
\frac{{\partial V\left( {{\bt^{0}},\bm} \right)}}{{\partial {t^{0}_a}}} &= \int_0^{x_a^{\text{user}}} {f\left( {\frac{s}{{{m_a}}}} \right)ds} , \label{cdc16-s1} \\
\frac{{\partial V\left( {{\bt^{0}},\bm} \right)}}{{\partial {m_a}}} &= \int_0^{x_a^{\text{user}}} {{t^{0}_a}\dot f\left( {\frac{s}{{{m_a}}}} \right)\left( { - \frac{s}{{m_a^2}}} \right)ds}, \label{cdc16-s2}
\end{align}
where $\dot f(\cdot)$ denotes the derivative of $f(\cdot)$. 

Note that typically we have $\frac{{\partial V\left( {{\bt^{0}},\bm} \right)}}{{\partial {t^{0}_a}}} > 0$ and $\frac{{\partial V\left( {{\bt^{0}},\bm} \right)}}{{\partial {m_a}}} < 0$, meaning
a slight decrease (resp., increase) of $t_a^0$ (resp., $m_a$) would reduce the objective function value of \eqref{cdc16-tap}. Based on this observation, for $a' \in \scrA$ we define the following quantities:
{\begin{align}
	& \Delta V\left( {{\bt^0},\bm;\Delta t_{a'}^0} \right) \defeq \mathop {\min }\limits_{\bx \in \scrF} \sum\limits_{a \in \scrA} {\int_0^{{x_a}} {t_a^0f\left( {\frac{s}{{{m_a}}}} \right)ds} }  \notag \\ 
	&- \mathop {\min }\limits_{\bx \in \scrF} \Bigg[ {\sum\limits_{a \in \scrA,a \ne a'} {\int_0^{{x_a}} {t_a^0f\left( {\frac{s}{{{m_a}}}} \right)ds} }}  \notag \\
	&~~~~~~~~~~+ \int_0^{{x_{a'}}} {\left( {t_{a'}^0 + \Delta t_{a'}^0} \right)f\left( {\frac{s}{{{m_{a'}}}}} \right)ds} \Bigg],   \label{sensi-freetime}
	\end{align}}
and
\begin{align}
	& \Delta V\left( {{\bt^0},\bm;\Delta m_{a'}} \right) \defeq \mathop {\min }\limits_{\bx \in \scrF} \sum\limits_{a \in \scrA} {\int_0^{{x_a}} {t_a^0f\left( {\frac{s}{{{m_a}}}} \right)ds} }  \notag \\ 
	&- \mathop {\min }\limits_{\bx \in \scrF} \Bigg[ {\sum\limits_{a \in \scrA,a \ne a'} {\int_0^{{x_a}} {t_a^0f\left( {\frac{s}{{{m_a}}}} \right)ds} }}  \notag \\
	&~~~~~~~~~~+ \int_0^{{x_{a'}}} {{t_{a'}^0} f\left( {\frac{s}{{{m_{a'} + \Delta m_{a'}}}}} \right)ds} \Bigg],   \label{sensi-capacity}
\end{align}
where {$\Delta t_{a'}^0 \defeq  - 0.2 \times \min \left\{ {t_a^0;\,a \in \scrA} \right\}$} and
{$\Delta {m_{a'}} \defeq 0.2 \times \min \left\{ {{m_a};\,a \in \scrA} \right\}$}. Note that, by construction, for each and every $a \in \scrA$, we approximately have $0 < \Delta V\left( {{\bt^0},\bm;\Delta t_{a}^0} \right) \propto {\partial V\left( {{\bt^{0}},\bm} \right)} / {{\partial {t^{0}_{a}}}}$ and $0 < \Delta V\left( {{\bt^0},\bm;\Delta m_{a}} \right) \propto \left| {\partial V\left( {{\bt^{0}},\bm} \right)} / {{\partial {m_a}}} \right|$.

\section{Dataset Description and Processing} 
\label{sec: dataset} 

In this section, based on our data-driven approach outlined in
Sec.~\ref{cdc16-sec:mod}, we conduct a case study using actual traffic
data from the Eastern Massachusetts (EMA) road network
\cite{BarGera16,InverseVIsTraffic}.

\subsection{Description of the Eastern Massachusetts dataset} \label{dataEMA}

\begin{figure}[htb]
	\centering
	\includegraphics[width=0.75\textwidth]{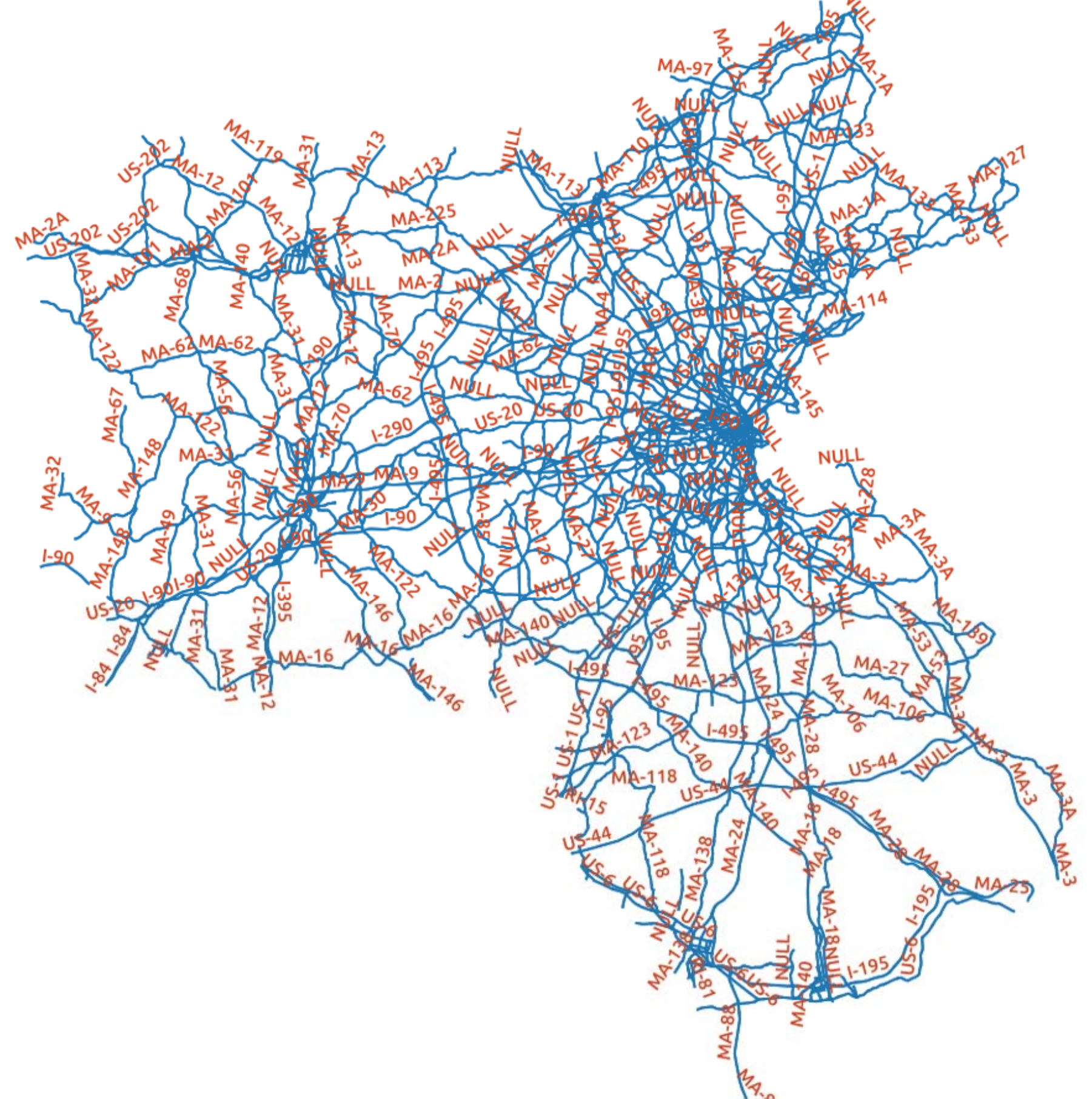}
	\caption{All available road segments in Eastern Massachusetts
          (from \cite{CDC16}).} 
	\label{eastMA}
\end{figure}

We deal with two datasets concerning the EMA road network: (\rmnum{1})
The speed dataset, made available to us by the Boston Region
Metropolitan Planning Organization (MPO), includes the spatial average
speeds for more than 13,000 road segments (with an average length of 0.7
miles; see Fig. \ref{eastMA}) of EMA, providing the average speed for
every minute of the year 2012.  For each road segment, identified with a
unique {\em tmc (traffic message channel)} code, the dataset provides
information such as speed data (instantaneous, average and free-flow
speed) in \emph{miles per hour (mph)}, date and time, and traveling time
(in \emph{minutes}) through that segment.
\noindent (\rmnum{2})
The flow capacity (in \emph{vehicles per hour}) dataset, also provided by the MPO, includes
capacity data for more than 100,000
road segments (with an average length of 0.13 miles) in EMA. For more detailed information of these two datasets, see \cite{CDC16}.

\begin{figure}[h]  
	\centering
	\begin{subfigure}[b]{0.4\textwidth}
		\includegraphics[width=\textwidth]{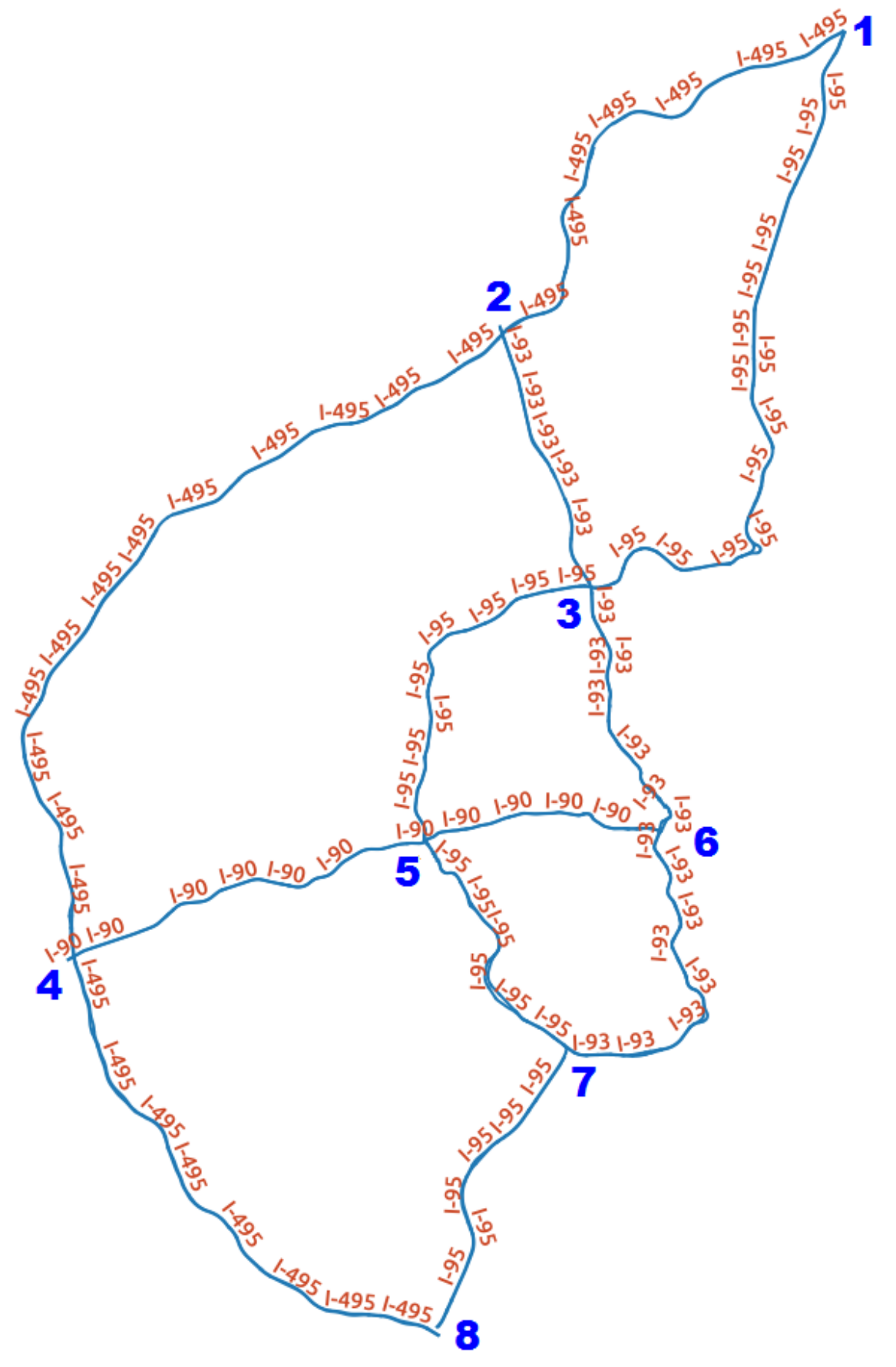}
		\caption{}
		\label{cdc16-sub-net-a}
	\end{subfigure} 
\begin{subfigure}[b]{0.4\textwidth}
		\includegraphics[width=\textwidth]{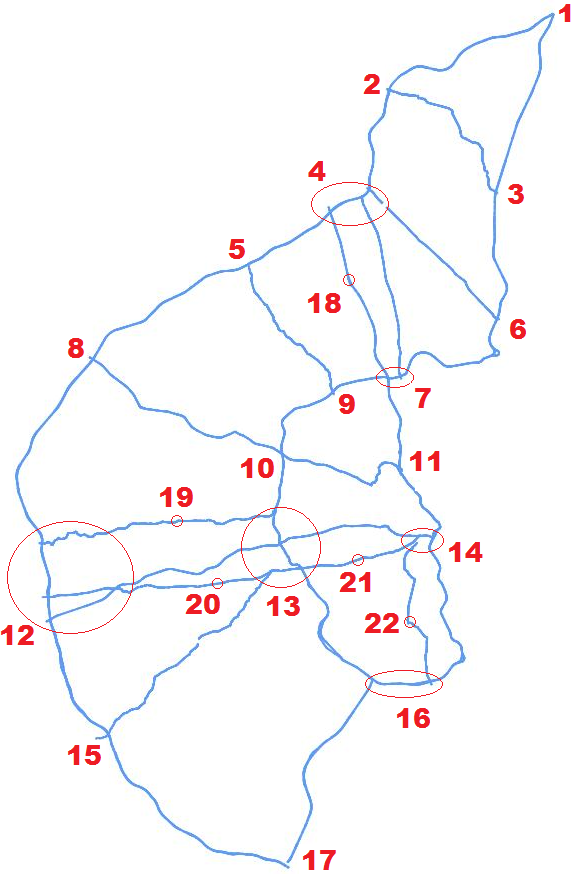}
		\caption{}
		\label{sub-net-a-ext}
	\end{subfigure}  
\caption{(a) An interstate highway sub-network of EMA ($\scrI_1$) (the blue numbers indicate node indices); (b) An extended highway sub-network of EMA ($\scrI_2$) (the red numbers indicate node indices). (See \cite{InverseVIsTraffic} for the correspondences between nodes and link indices.)}
	\label{cdc16-sub-net}
\end{figure}

\begin{figure}[htp]  
	\centering
	\includegraphics[width=0.95\textwidth]{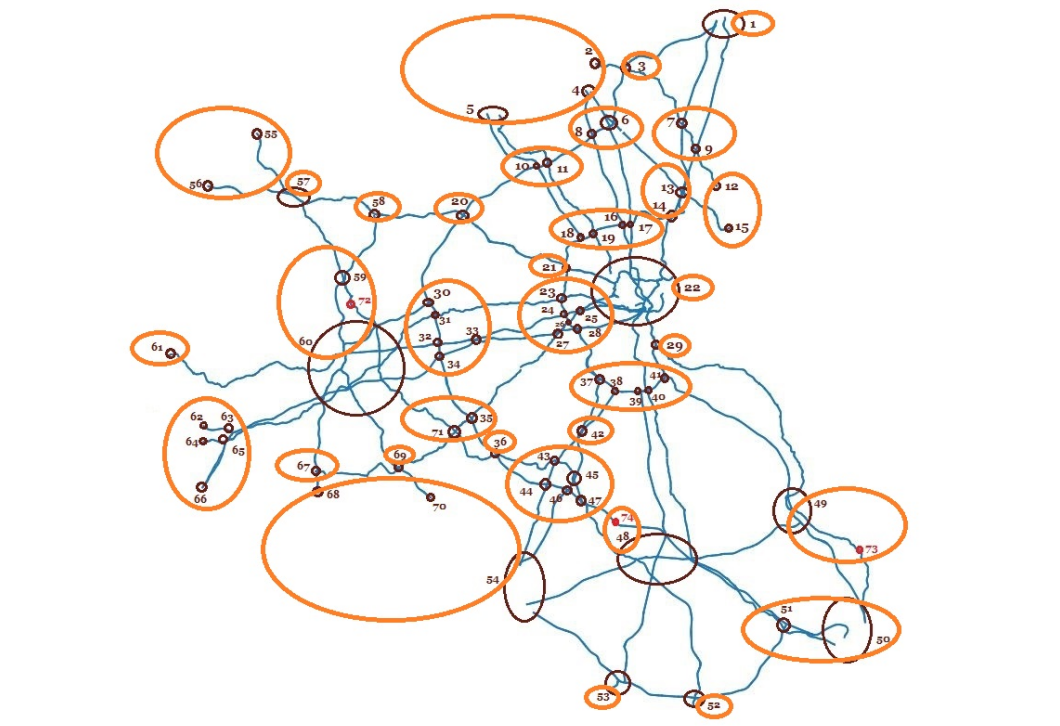}
\caption{A wider EMA highway subnetwork ($\scrI_3$); details on the
  correspondences between nodes and link indices are in
  \cite{InverseVIsTraffic}. (``nodes:zone'' pairs -- \{1\}: Seabrook (NH);
  \{2, 4, 5\}: NH; \{3\}: Haverhill; \{6, 8\}: Lawrence; \{7, 9\}: Georgetown;
  \{10, 11\}: Lowell; \{12, 15\}: Salem; \{13, 14\}: Peabody; \{16, 17, 18,
  19\}: Burlington; \{20\}: Littleton; \{21\}: Lexington; \{22\}: Boston;
  \{23, 24, 25, 26, 27, 28\}: Waltham; \{29\}: Quincy; \{30, 31, 32, 33, 34\}:
  Marlborough/Framingham; \{35, 71\}: Milford; \{36\}: Franklin; \{37, 38, 39,
  40, 41\}: Westwood/Quincy; \{42\}: Dedham; \{43, 44, 45, 46, 47\}:
  Foxborough; \{48, 74\}: Taunton; \{49, 73\}: Plymouth; \{50, 51\}: Cape Cod;
  \{52\}: Dartmouth; \{53\}: Fall River; \{54, 68, 70\}: RI; \{55, 56\}: VT;
  \{57\}: Westminster; \{58\}: Leominster; \{59, 60, 72\}: Worcester; \{61\}:
  Amherst; \{62, 63, 64, 65, 66\}: CT; \{67\}: Webster; \{69\}: Uxbridge.)}
	\label{fig:zoneMA}
\end{figure}

\subsection{Preprocessing} \label{prep}

In \cite{CDC16} and \cite{poa-ifac-2017} we investigate two relatively small
subnetworks (denoted by $\scrI_1$ and $\scrI_2$ and shown in
Figs. \ref{cdc16-sub-net-a} and \ref{sub-net-a-ext}, respectively) of
the EMA road network. Here, we further consider a much larger subnetwork
(denoted by $\scrI_3$ and shown in Fig.~\ref{fig:zoneMA}).  Performing
similar preprocessing procedures as those in \cite{CDC16,poa-ifac-2017}, we end
up with traffic flow data (Wardrop equilibria) and road (link)
parameters (\textit{flow capacity} and \textit{free-flow travel time})
for the three subnetworks $\scrI_1$, $\scrI_2$, and $\scrI_3$, where
$\scrI_1$ contains only interstate highways, $\scrI_2$ also contains
state highways, and $\scrI_3$ covers a much wider area of EMA. Note that
$\scrI_1$ (resp., $\scrI_2$, $\scrI_3$) consists of 8 (resp., 22, 74)
nodes and 24 (resp., 74, 258) links. Assuming that each node could be an
origin and a destination, then there are $8 \times (8-1) = 56$ (resp.,
$22 \times (22-1) = 462$) OD pairs in $\scrI_1$ (resp., $\scrI_2$). For
$\scrI_3$, we simplify the analysis by grouping nodes within the same
area, assigning them the same \emph{zone} label, thus obtaining 34 zones
(as opposed to 74 nodes). Assuming that each zone could be an origin and
a destination, then there are $34 \times (34-1) = 1122$ OD pairs in
$\scrI_3$. It is worth pointing out that nodes 18, 19, 20, 21, and 22
(resp., 72, 73, and 74) in $\scrI_2$ (resp., $\scrI_3$) are introduced
for ensuring the identifiability of the OD demand matrices. More
specifically, to ``recover'' uniquely an OD demand matrix from observed
link flow data, the link-route incidence matrix $\bA$ is required to
satisfy certain structural properties; see \cite[Lemma
2]{hazelton2000estimation}.

\subsection{Estimating initial OD demand matrices} \label{sec:ODmat}

Operating on $\scrI_1$, we solve the QP (P1) (cf. \eqref{qp2}) and the
QCP (P2) (cf. \eqref{qp3}) using data corresponding to five different
time periods (AM, MD (middle day), PM, NT (night), and WD (weekend)) of
four months (Jan., Apr., Jul., and Oct.) in 2012, thus obtaining 20
different OD demand matrices for these scenarios. Expanding each and
every of the 20 OD demand matrices of $\scrI_1$ by setting the demand
for any OD pair that belongs to $\scrI_2$ but does not belong to
$\scrI_1$ to zero, we obtain ``rough'' initial demand matrices for
$\scrI_2$.

On the other hand, for the much larger subnetwork $\scrI_3$, to obtain
initial OD demand matrices corresponding to the same 20 scenarios, we
perform a different simplification procedure. In particular, we only
consider the \emph{shortest} route for each OD pair of $\scrI_3$,
thus leading 
to a deterministic route choice matrix $\bP$ and significantly reducing
the number of decision variables in the QCP (P2).

As noted in \cite{CDC16}, the GLS method assumes the traffic network to
be uncongested. It follows, that the estimated OD demand matrices for
non-peak periods (MD/NT/WD) are relatively more accurate than those for
peak periods (AM/PM). After obtaining estimates for travel latency cost
functions in Sec.~\ref{cost}, based on the observed Wardrop flows and the initial estimates for the OD demand matrices, we will conduct
demand adjustment procedures for $\scrI_2$ and $\scrI_3$ in
Sec.~\ref{od-adj}.

\subsection{{Cost function estimation and sensitivity analysis}} \label{cost}

First, to validate the effectiveness and efficiency of the cost function
estimator \eqref{cdc16-costEstimator}, we conduct numerical experiments
over the Anaheim benchmark network \cite{BarGera16}, whose ground truth
cost functions, OD demand matrices, and all necessary road parameters
are available.  Next, operating on $\scrI_1$ using the flow data and the
OD demand matrices obtained in Secs.~\ref{prep} and \ref{sec:ODmat}
respectively, we estimate the travel latency cost functions, $f(\cdot)$
in particular, for 20 different scenarios, via the estimator
\eqref{cdc16-costEstimator}, by solving the QP invVI-2 accordingly. As
in \cite{CDC16}, to make the estimates reliable, for each scenario, we
perform a 3-fold cross-validation.  Note that \cite{CDC16} applied a
different estimator, which is numerically not as stable.

We assume that such estimates for $f(\cdot)$, as obtained from $\scrI_1$,
can be shared by all the three subnetworks $\scrI_1$, $\scrI_2$, and
$\scrI_3$; this makes sense, because the function $f(\cdot)$ is common
for all links and, when estimating it through $\scrI_1$, we have already
made use of a large amount of data (note that there are 24 links in
$\scrI_1$ and the flow data and the corresponding OD demand matrices
that we use have covered $120$ different time instances for each of the
20 scenarios; for details, see \cite{InverseVIsTraffic}).

To illustrate our method of analyzing sensitivities for the TAP
formulation \eqref{cdc16-tap}, we again conduct numerical experiments on
$\scrI_1$. In particular, we investigate a scenario corresponding to the
AM peak period of April 2012.

\subsection{OD demand adjustments} \label{od-adj}

First, we demonstrate the effectiveness of
Alg.~\ref{alg:demandAdjustment} using the Anaheim benchmark network.
Then, assuming the per-road travel latency cost functions are available
(we take the travel latency cost functions derived from $\scrI_1$ as in
Sec. \ref{cost}), we apply Alg. \ref{alg:demandAdjustment} to $\scrI_2$,
which contains $\scrI_1$ as one of its representative subnetworks. Note
that the main difference between $\scrI_1$ and $\scrI_2$ is the modeling
emphasis; specifically, $\scrI_1$ only takes account of interstate
highways, while $\scrI_2$ also encompasses state highways, thus
containing more details of the real road network of EMA. We can think of
$\scrI_1$ as a ``landmark'' subnetwork of $\scrI_2$. Based on the
initially estimated demand matrices for $\scrI_1$, we will implement the
following generic demand-adjusting scheme so as to derive the OD demand
matrices for $\scrI_2$.

Given a network ($\scrI_2$ in our case) of any size we can select its
``landmark'' subnetworks ($\scrI_1$ in our case) (based on the
information of road types, pre-identified centroids, etc.) with
acceptably smaller sizes; say we end up with $N$ ($N=1$ in our case)
such subnetworks. Then, for each subnetwork, we estimate its demand
matrix by solving sequentially the QP (P1) and the QCP (P2)
(cf. Sec. \ref{sec:ODmat}). Setting the demand for any OD pair not
belonging to this subnetwork to zero, we obtain a ``rough'' initial
demand matrix for the entire network ($\scrI_2$ in our case). Next, we take
the average of these initial demand matrices. Finally, we adjust the
average demand matrix based on the flow observations of the entire
network.

Next, taking again the travel latency cost functions derived from
$\scrI_1$, we apply Alg. \ref{alg:demandAdjustment} to $\scrI_3$, based
on the initial OD demand matrices estimated from $\scrI_3$ (see
Sec. \ref{sec:ODmat}) and the Wardrop flows inferred from $\scrI_3$ (see
Sec. \ref{prep}).

As noted in Remark \ref{rem:od}, the reason for not directly solving
(P1) and (P2) for the larger networks ($\scrI_2$ and $\scrI_3$ in our
case) is that there are too many decision variables in (P2) and this
would lead to numerical difficulties.

\subsection{{PoA evaluation and meta analysis}}

We calculate the PoA values for $\scrI_2$ and $\scrI_3$ for the PM
period of April 2012. For each day, in \eqref{poa-def} we take the
average observed link flows over the PM period as the ``user flows,''
and obtain the ``social flows'' by solving the NLP \eqref{obj_soc} using
the estimated cost functions and demand matrix exclusively for the PM
period. To solve \eqref{obj_soc}, we use the IPOPT solver \cite{ipopt}
which implements a primal-dual interior point method
\cite{wright1997primal}.

To better understand the performance of the road network under the
user-centric vs. the system-centric routing policy, we conduct a meta
analysis on $\scrI_3$. In particular, under the two policies, we compare
congestion for various zones of the network, the maximum/minimum link
flows, and link-specific congestion.

\section{Numerical Results} \label{Sec:Rsults}

For economy of space, we will not show the detailed results for the
initial estimation of OD demand matrices. {However, we {report in
    Tab.~\ref{journal-tab2}} the entries of the route choice probability
  matrix $\bP$ derived for $\scrI_1$ for some specific OD pairs (the
  complete results can be found in \cite{InverseVIsTraffic}). {It is
    seen from Tab.~\ref{journal-tab2} that we cannot always expect a
    higher probability for a shorter/faster route; randomness
    exists. However, this is not necessarily counterintuitive, because
    the selected routes for the same OD pair have close lengths/travel
    times.}}  {We note here that {when identifying (and refining) the
    feasible routes for each OD pair of $\scrI_1$, we consider at most
    three \emph{shortest} routes (ranked \#1-\#3) and discard all the
    routes with a length larger than that of the \emph{shortest} route
    (ranked \#1) by more than $50\%$. Note also that since this initial
    OD estimation procedure does not involve real-time updates of
    traffic conditions, we may use either travel times or lengths as
    weights for links in the graph model.}}

\begin{table}[t]{
	\newcommand{\tabincell}[2]{\begin{tabular}{@{}#1@{}}#2\end{tabular}}
	\centering
	\caption{Selected route choice analysis results for
		$\scrI_1$ (corresponding to the PM peak period of April 2012).} \label{journal-tab2} 
	\resizebox{0.98\textwidth}{!}
	{
		\begin{tabular}{clccc}
			\midrule[1.5pt]
			\textbf{OD pair}  & \textbf{refined feasible route}  & \tabincell{c}{\textbf{route length} \\\textbf{(in miles)}} & \tabincell{c}{\textbf{free-flow travel time} \\\textbf{(in hours)}} & \tabincell{l}{\textbf{route choice}\\ \textbf{probability}}  \\
			\midrule[1.5pt]
			(1, 8)  & \tabincell{l}{$1 \to 2 \to 3 \to 5 \to 7 \to 8$ \\ $1 \to 3 \to 5 \to 7 \to 8$\\ $1 \to 2 \to 3 \to 6 \to 7 \to 8$} & \tabincell{l}{74.2072 \\ 74.6696 \\ 74.8692} & \tabincell{l}{1.0235 \\ 1.0297 \\ 1.0522}  &   \tabincell{l}{0.3265 \\ 0.3394 \\ 0.3341}  \\
			\midrule
			(2, 4)  & \tabincell{l}{$2 \to 4$ \\ $2 \to 3 \to 5 \to 4$ \\ $2 \to 3 \to 6 \to 5 \to 4$} & \tabincell{l}{37.6346 \\ 43.4554 \\ 50.7995} & \tabincell{l}{0.5123 \\ 0.6010 \\ 0.7274}  &  \tabincell{l}{0.8274 \\ 0.1004 \\ 0.0722} \\
			\midrule
			(3, 5)  & \tabincell{l}{$3 \to 5$ \\$3 \to 6 \to 5$} & \tabincell{l}{16.2154 \\ 23.5596} & \tabincell{l}{0.2262 \\ 0.3526}   &  \tabincell{l}{0.8375 \\ 0.1625}\\
			\midrule
			(8, 3)  & \tabincell{l}{$8 \to 7 \to 5 \to 3$ \\$8 \to 7 \to 6 \to 3$ \\ $8 \to 7 \to 5 \to 6 \to 3$} & \tabincell{l}{43.3260 \\ 43.4313 \\50.2382} & \tabincell{l}{0.6065 \\ 0.6310 \\ 0.7308}   &  \tabincell{l}{0.4364 \\ 0.3022 \\ 0.2614}\\
			\midrule[1.5pt]
		\end{tabular}}}
	\end{table}

In the following, we will focus on presenting the results for the
estimates of the travel latency cost functions (derived for the Anaheim
benchmark network and $\scrI_1$), the demand adjustment procedure
(derived for the Anaheim benchmark network; note that we will not show
the detailed demand adjustment results for $\scrI_2$ and $\scrI_3$,
because we do not have the ground truth for a comparison), the PoA
evaluations (derived for $\scrI_2$ and $\scrI_3$), the sensitivity
analysis (derived for $\scrI_1$), and the meta analysis (derived for
$\scrI_3$).

\subsection{Results from estimating the cost functions}

\subsubsection{Results for the Anaheim benchmark network}

\begin{figure}[h]  
	\centering
	\begin{subfigure}[b]{0.49\textwidth}
		\includegraphics[width=\textwidth]{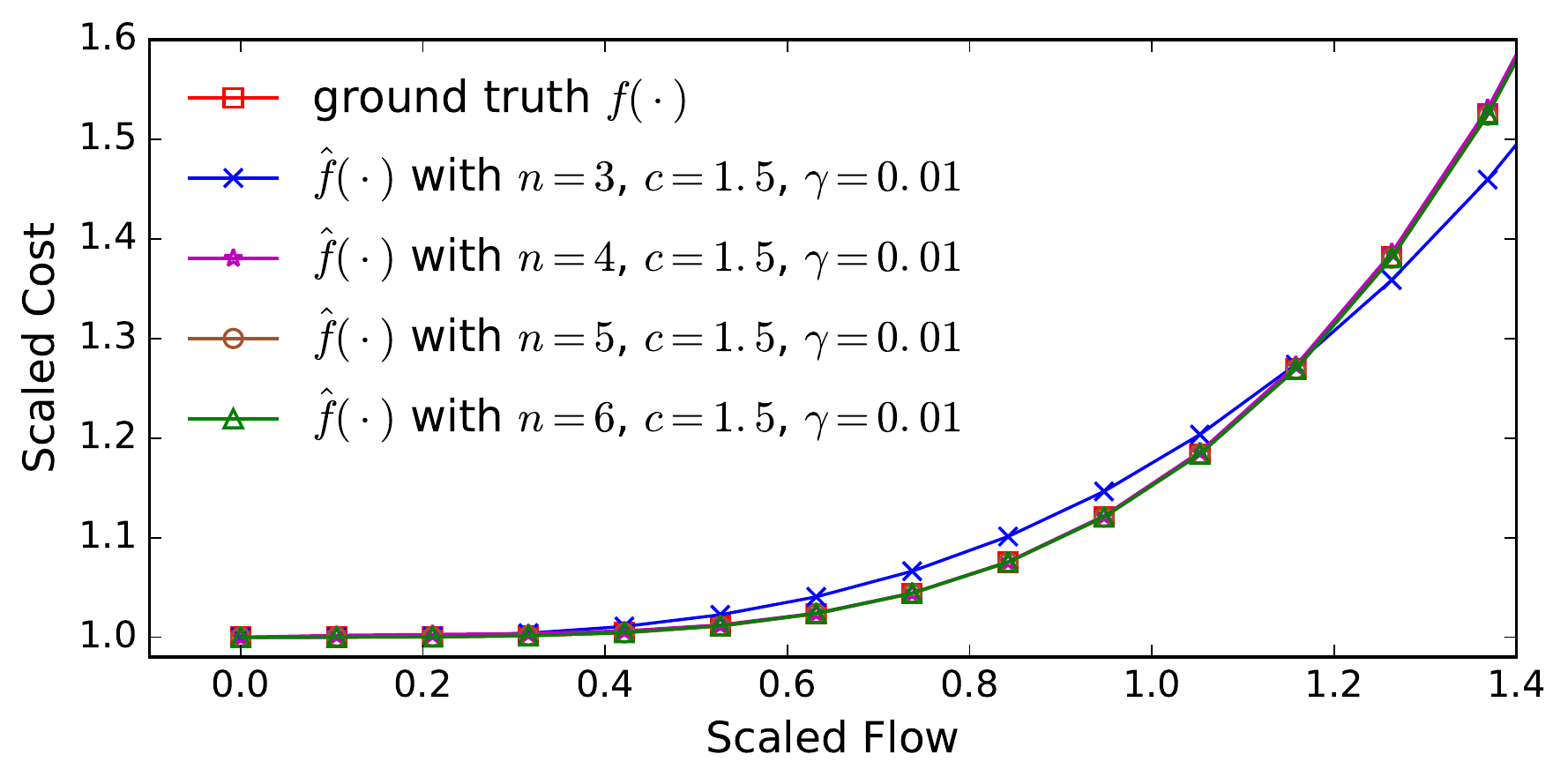}
		\caption{Vary $n$ ($c$ and $\gamma$ fixed).}
		\label{fig:n_Anaheim}
	\end{subfigure} 
	\begin{subfigure}[b]{0.49\textwidth}
		\includegraphics[width=\textwidth]{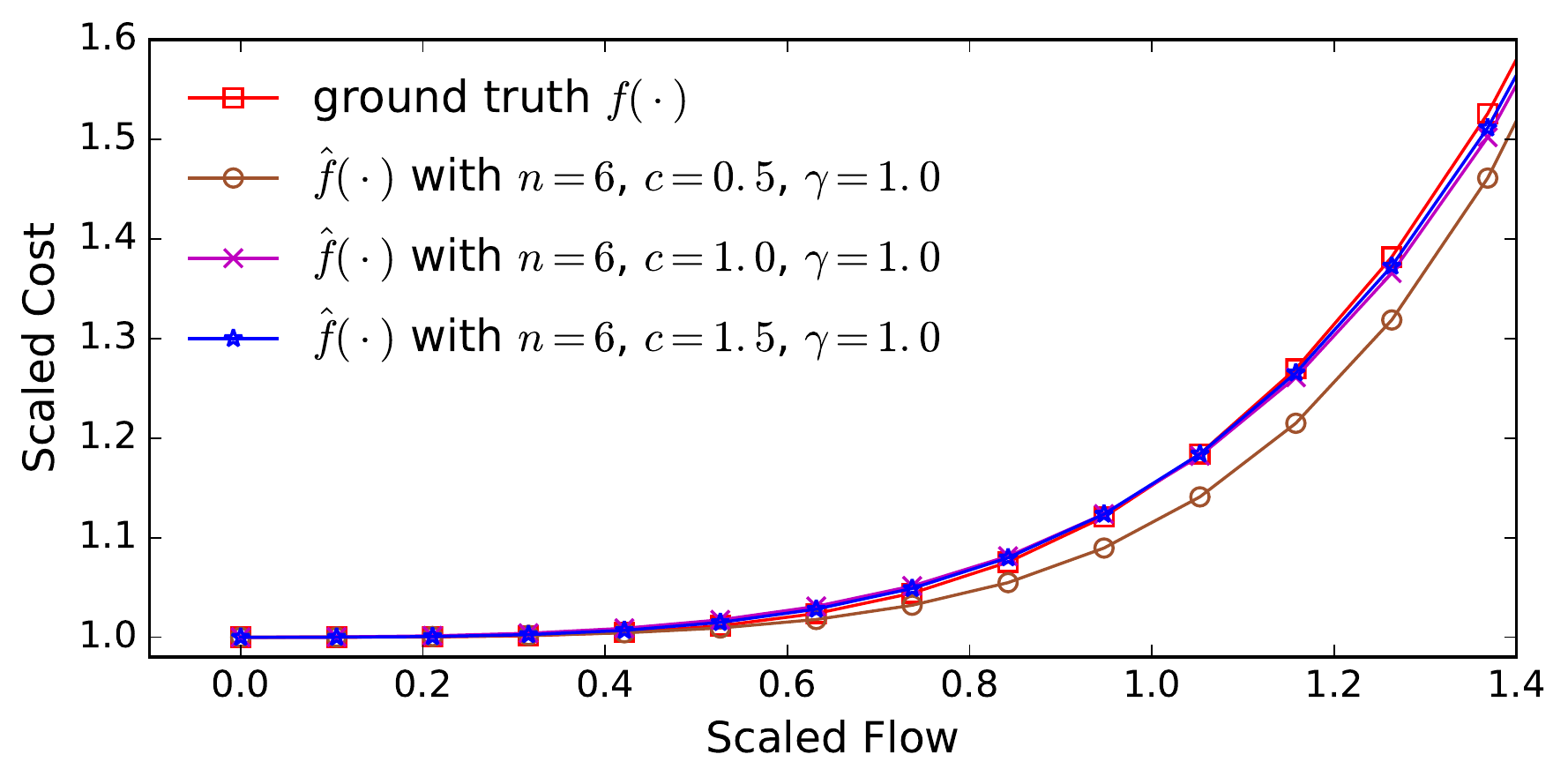}
		\caption{Vary $c$ ($n$ and $\gamma$ fixed).}
		\label{fig:c_Anaheim}
	\end{subfigure}   
	\begin{subfigure}[b]{0.49\textwidth}
		\includegraphics[width=\textwidth]{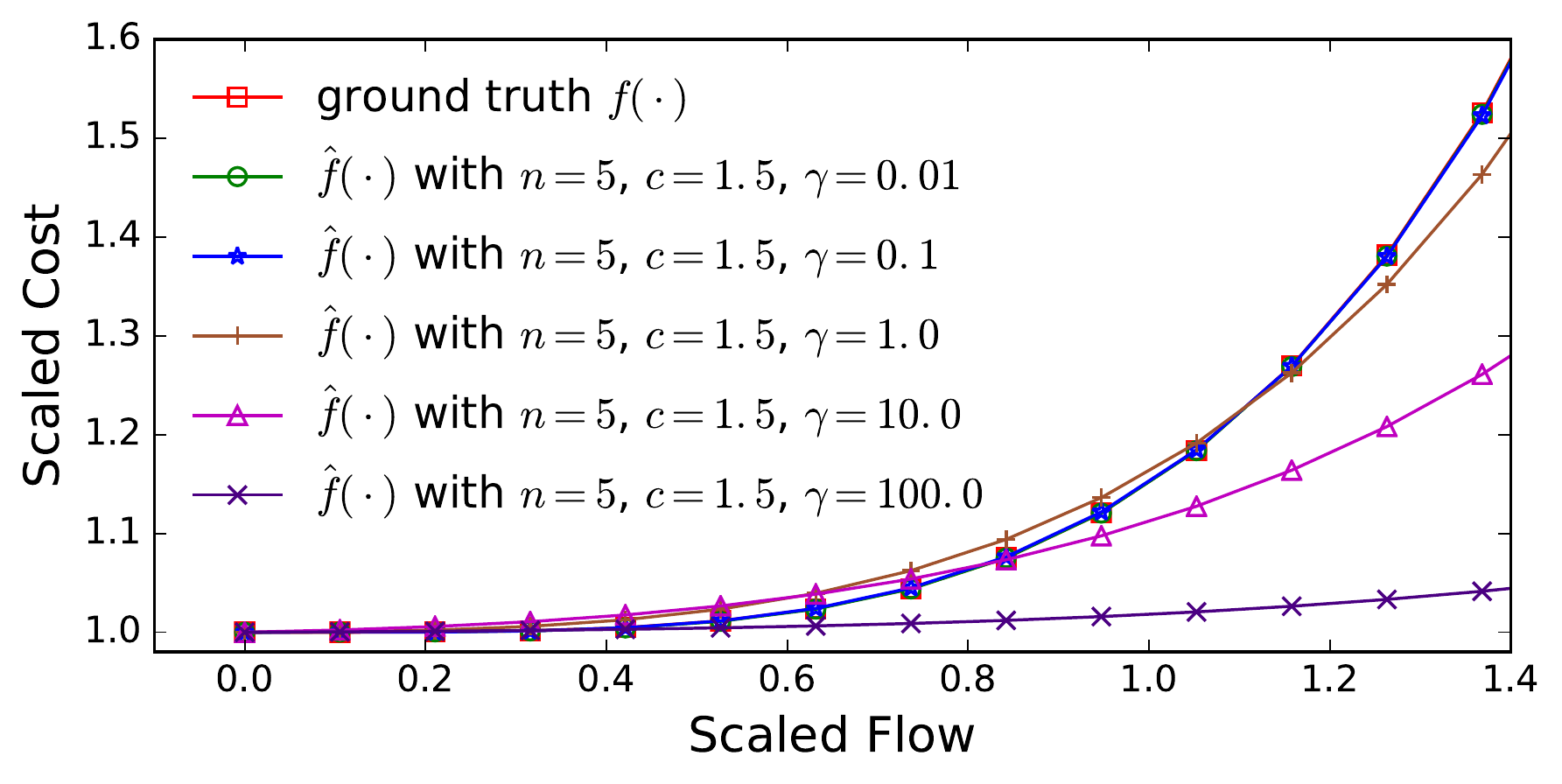}
		\caption{Vary $\gamma$ ($n$ and $c$ fixed).}
		\label{fig:gamma_Anaheim}
	\end{subfigure}  
	\caption{Estimations for cost function $f(\cdot)$ by solving invVI-2 corresponding to different parameter settings (Anaheim).}
	\label{fig:cost_Anaheim}
\end{figure}

The Anaheim network contains 38 zones (hence $38 \times (38 - 1) = 1406$
OD pairs), 416 nodes, and 914 links. The ground truth
$f\left(\cdot\right)$ is taken as $f\left(z\right) = 1 + 0.15z^4, \,z
\geq 0$. Fig. \ref{fig:cost_Anaheim} shows the estimation results for
$f(z)$ by solving invVI-2 corresponding to different parameter
settings. In particular, Fig. \ref{fig:n_Anaheim} shows the curves of
the ground truth $f(z)$ and the estimator $\hat f(z)$ corresponding to
$n$ taking values from $\{3, 4, 5, 6\}$ while keeping $c$ and $\gamma$
fixed to 1.5 and 0.01 respectively; it is seen that except for the case
$n = 3$, all estimation curves are very close to the ground truth. Note
that the ground truth $f(z)$ is a polynomial function with degree 4,
which is greater than 3. This suggests the use of a value $n\geq 4$ in
recovering the cost function $f(\cdot)$. The intuition here is that we
can use a higher order polynomial with appropriate coefficients to
approximate a lower order polynomial, but not vice versa.
Fig. \ref{fig:c_Anaheim} shows the curves of the ground truth $f(z)$ and
the estimator $\hat f(z)$ corresponding to $c$ taking values from
$\{0.5, 1.0, 1.5\}$ while keeping $n$ and $\gamma$ fixed to 6 and 1.0
respectively; it is seen that except for the case $c = 0.5$, the
estimation curves are very close to the ground truth. This suggests that
setting $c$ reasonably larger should give better estimation
results. Fig. \ref{fig:gamma_Anaheim} plots the curves of the ground
truth $f(z)$ and the estimator $\hat f(z)$ corresponding to $\gamma$
taking values from $\{0.01, 0.1, 1.0, 10.0, 100.0\}$ while keeping $n$
and $c$ fixed to 5 and 1.5 respectively; it is seen that as $\gamma$ is
set smaller and smaller, the estimation curve gets closer and closer to
the ground truth. This suggests that choosing a smaller regularization
parameter $\gamma$ should give tighter estimation results in terms of
data reconciling.

\subsubsection{Results for $\scrI_1$}

\begin{figure}[h]  
	\centering
	\begin{subfigure}[b]{0.49\textwidth}
		\includegraphics[width=\textwidth]{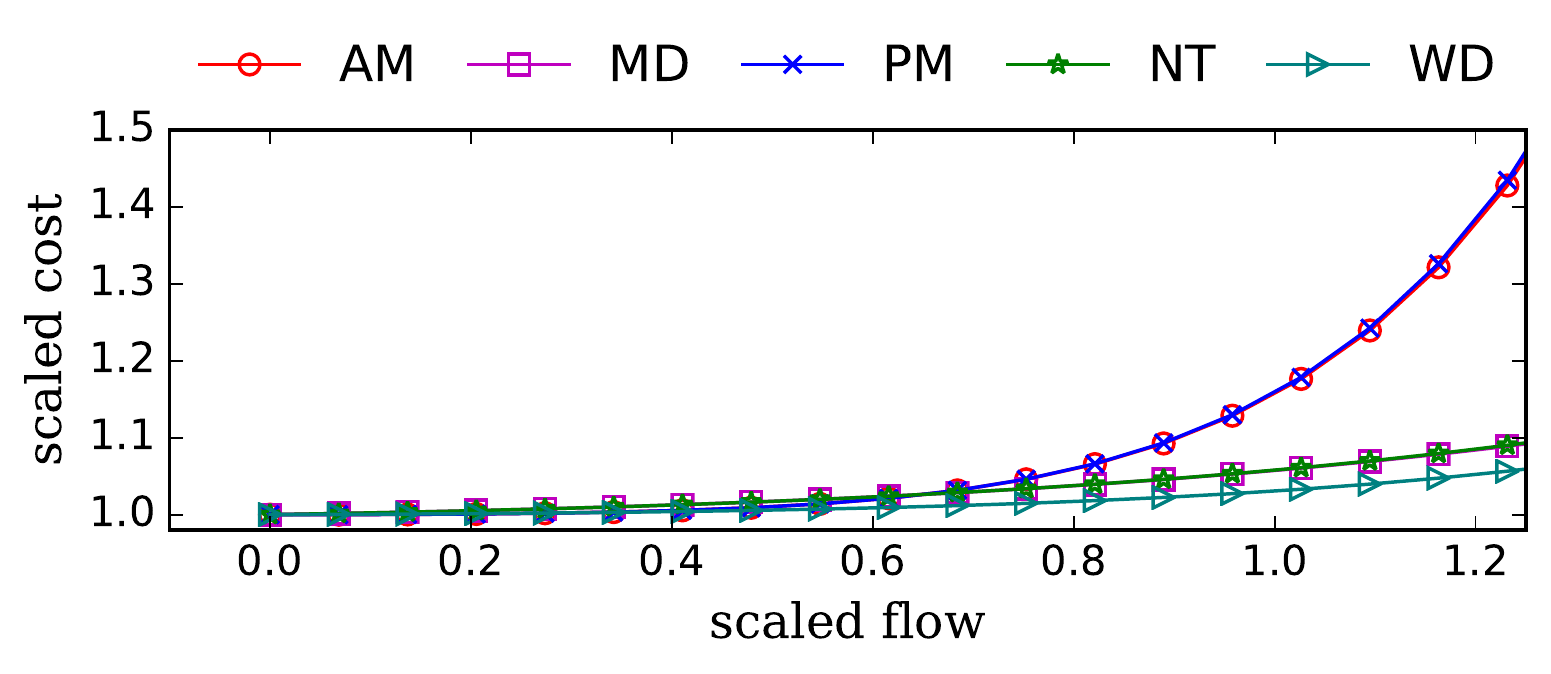}
		\caption{Jan.}
		\label{comp-cost-a}
	\end{subfigure} 
	\begin{subfigure}[b]{0.49\textwidth}
		\includegraphics[width=\textwidth]{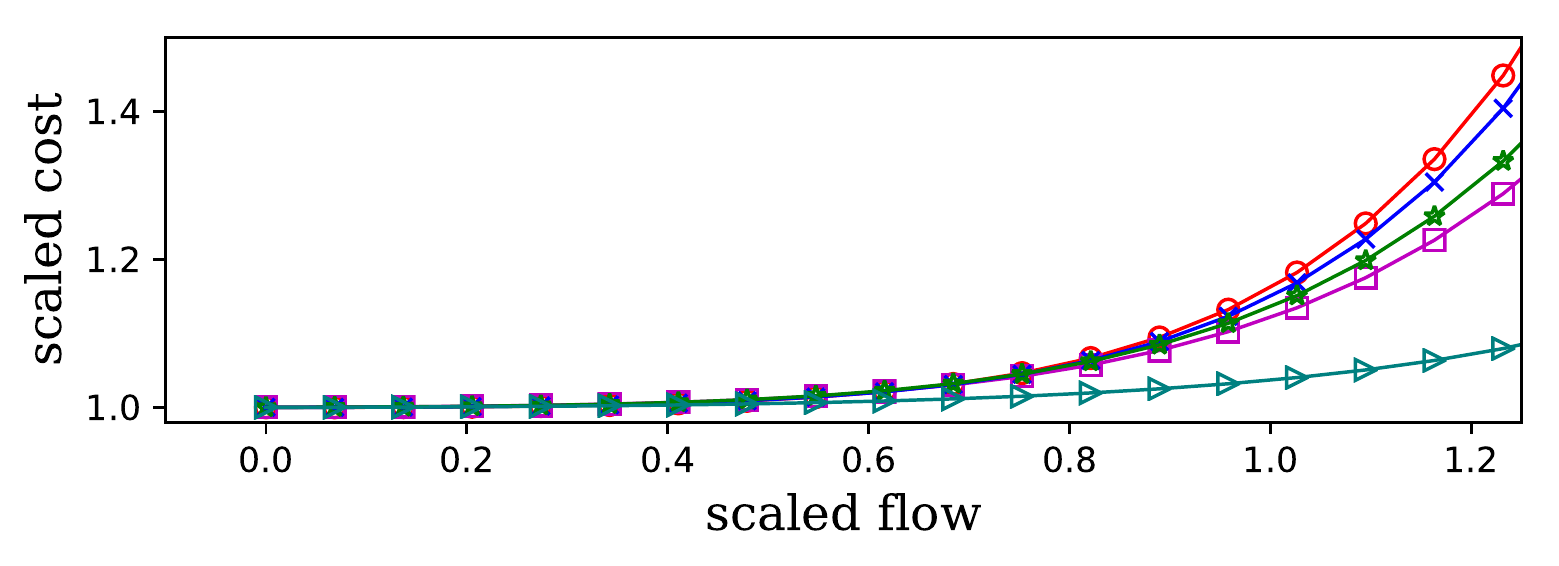}
		\caption{Apr.}
		\label{comp-cost-b}
	\end{subfigure}   
	\begin{subfigure}[b]{0.49\textwidth}
	\includegraphics[width=\textwidth]{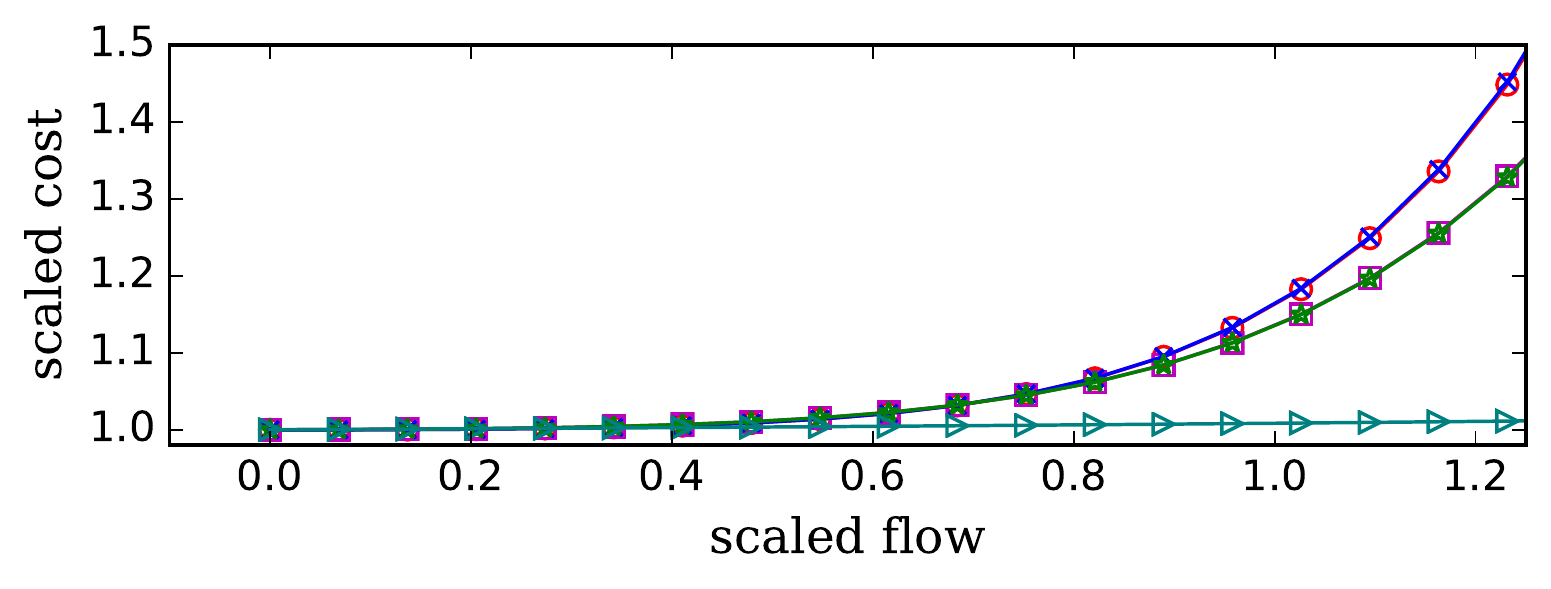}
	\caption{Jul.}
	\label{comp-cost-c}
	\end{subfigure}  
	\begin{subfigure}[b]{0.49\textwidth}
	\includegraphics[width=\textwidth]{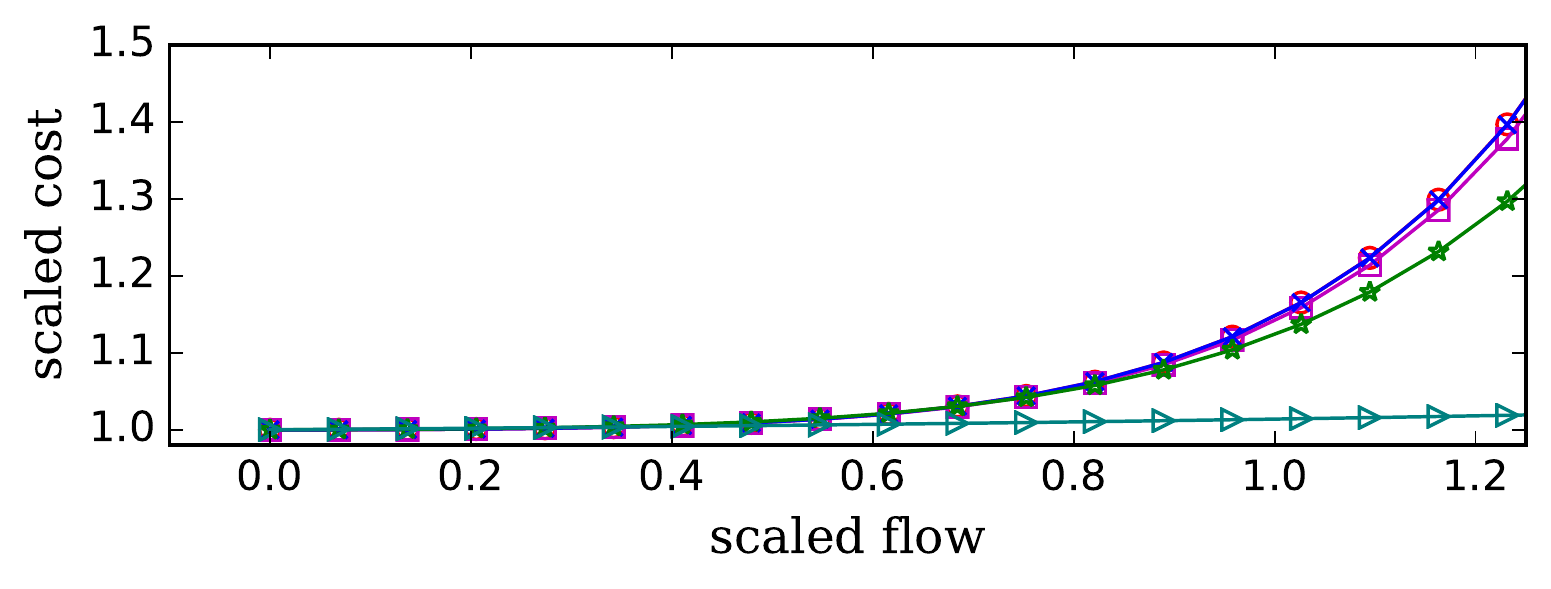}
	\caption{Oct.}
	\label{comp-cost-d}
	\end{subfigure}  
	\caption{Estimates for $f(\cdot)$ corresponding to different time periods (AM, MD (middle day), PM, NT (night), WD (weekend)), derived from data over $\scrI_1$ for 2012.}
	\label{comp-cost}
\end{figure}

We show the comparison results of the cost functions in
Fig. \ref{comp-cost}, where in each sub-figure, we plot the curves of
the estimated $f(\cdot)$ corresponding to five different time periods. For economy of space, we will not list the parameter setting details of $n$, $c$, and $\gamma$, which were selected by conducting a 3-fold cross-validation.

We observe from Figs.~\ref{comp-cost-a}-\ref{comp-cost-d} that the
costs for peak periods (AM/PM) are more sensitive to traffic flows than
for non-peak periods (MD/NT/WD). This can be
explained as follows: during rush hour, it is very common for vehicles
to pass through a congested road network while during non-rush hour, drivers mostly enjoy an uncongested road network.

In addition, it is seen that, for different months, the cost curves
for non-peak periods differ more significantly than for peak periods. Aside from the observation and
modeling errors, this can also be explained by seasonal traveling
patterns.

\subsection{Results from OD demand adjustment}  \label{sec:od-adj}

\begin{figure}[h]  
	\centering
	\begin{subfigure}[b]{0.49\textwidth}
		\includegraphics[width=\textwidth]{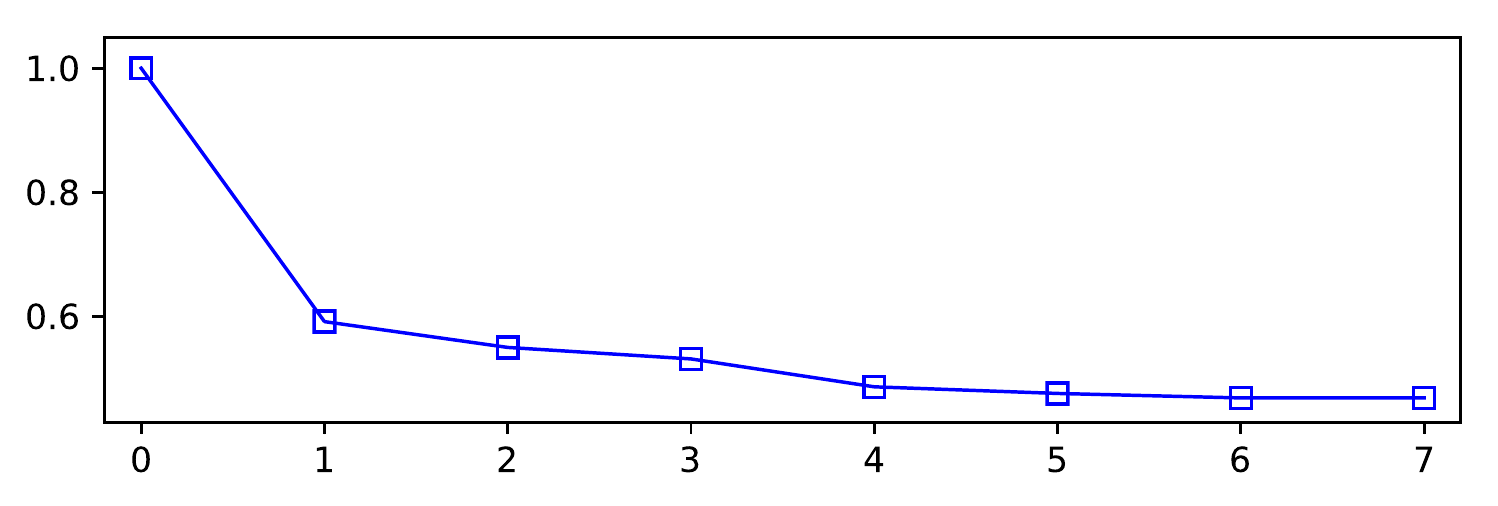}
		\caption{$F(\bg^l)/F(\bg^0)$ vs. \# of iterations.}
		\label{objFun_Anaheim}
	\end{subfigure} 
	\begin{subfigure}[b]{0.49\textwidth}
		\includegraphics[width=\textwidth]{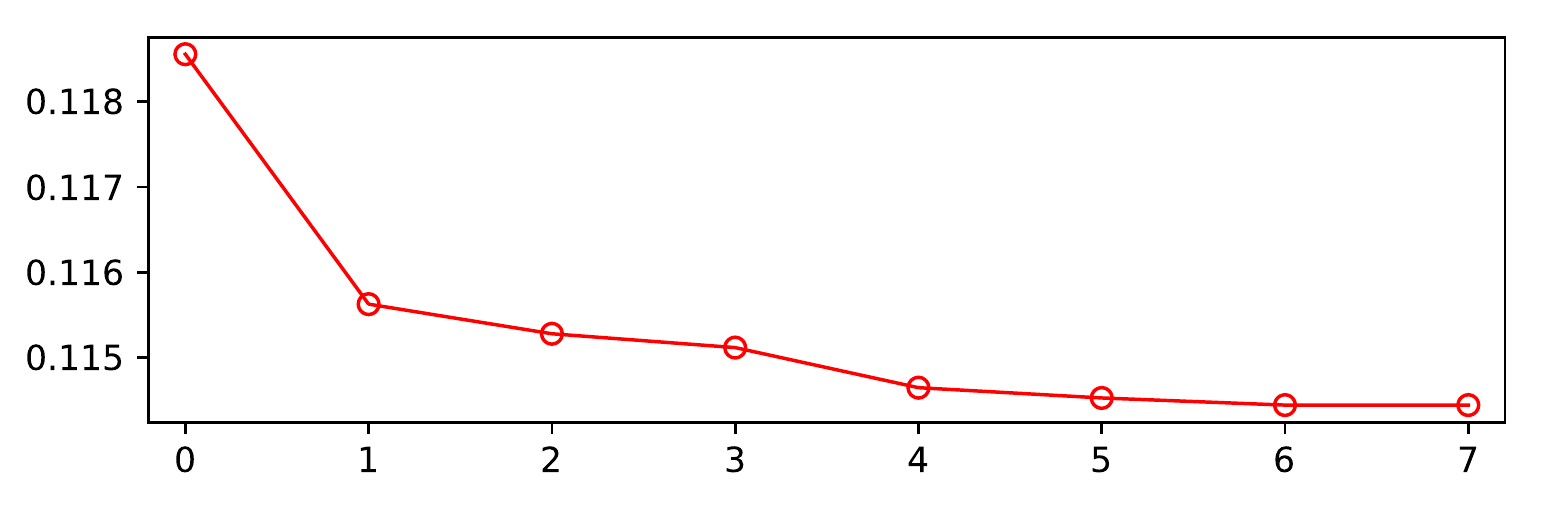}
		\caption{$\|\bg^l - \bg^*\| / \| \bg^* \|$ vs. \# of iterations.}
		\label{demandsDiff_biLev_Anaheim}
	\end{subfigure}   
	\caption{Key quantities vs. \# of iterations (Anaheim).}
	\label{fig:Anaheim}
\end{figure}

{We now present the OD demand adjustment results from the Anaheim
  network. For each OD pair, the initial demand is taken by scaling the
  ground truth demand using a random factor with uniform distribution over
  $[0.8, 1.2]$. The ground truth $f\left(\cdot\right)$ is taken as
  $f\left(z\right) = 1 + 0.15z^4, \; \forall z \geq 0$, and is assumed
  directly available.  When implementing
  Alg. \ref{alg:demandAdjustment}, we set $\gamma_1 = 0$, $\gamma_2 =
  1$, $\rho = 2$, {$T = 10$}, $\varepsilon_1 =0$, and $\varepsilon_2 =
  10^{-20}$. Fig.~\ref{objFun_Anaheim} shows that, after 7 iterations,
  the objective function value of the BiLev \eqref{bilevel} has been
  reduced by more than $50\%$. Fig.~\ref{demandsDiff_biLev_Anaheim}
  shows that, the distance between the adjusted demand and the ground
  truth demand keeps decreasing with the number of iterations, and the
  distance changes very slightly, meaning the adjustment procedure does
  not alter the initial demand much.}  Note that in
Fig. \ref{objFun_Anaheim}, the vertical axis corresponds to the
normalized objective function value of the BiLev, i.e.,
$F(\bg^l)/F(\bg^0)$ and, in Fig. \ref{demandsDiff_biLev_Anaheim}, the
vertical axis denotes the normalized distance between the adjusted
demand vector and the ground truth, i.e., $\|\bg^l - \bg^*\| / \| \bg^*
\|$, where $\bg^{*}$ is the ground-truth demand vector.

\subsection{Results for PoA evaluation}

After implementing the demand adjusting scheme, we obtain the demand
matrices for $\scrI_2$ and $\scrI_3$ on a daily-basis, as opposed to
those for $\scrI_1$ on a monthly-basis. Note that, even for the same
period of a day and within the same month, slight demand variations
among different days are possible; thus, our PoA results for $\scrI_2$
and $\scrI_3$ would be more accurate than those for $\scrI_1$ (shown in
\cite{CDC16}).
		
\begin{figure}[h]  
	\centering
	\begin{subfigure}[b]{0.49\textwidth}
		\includegraphics[width=\textwidth]{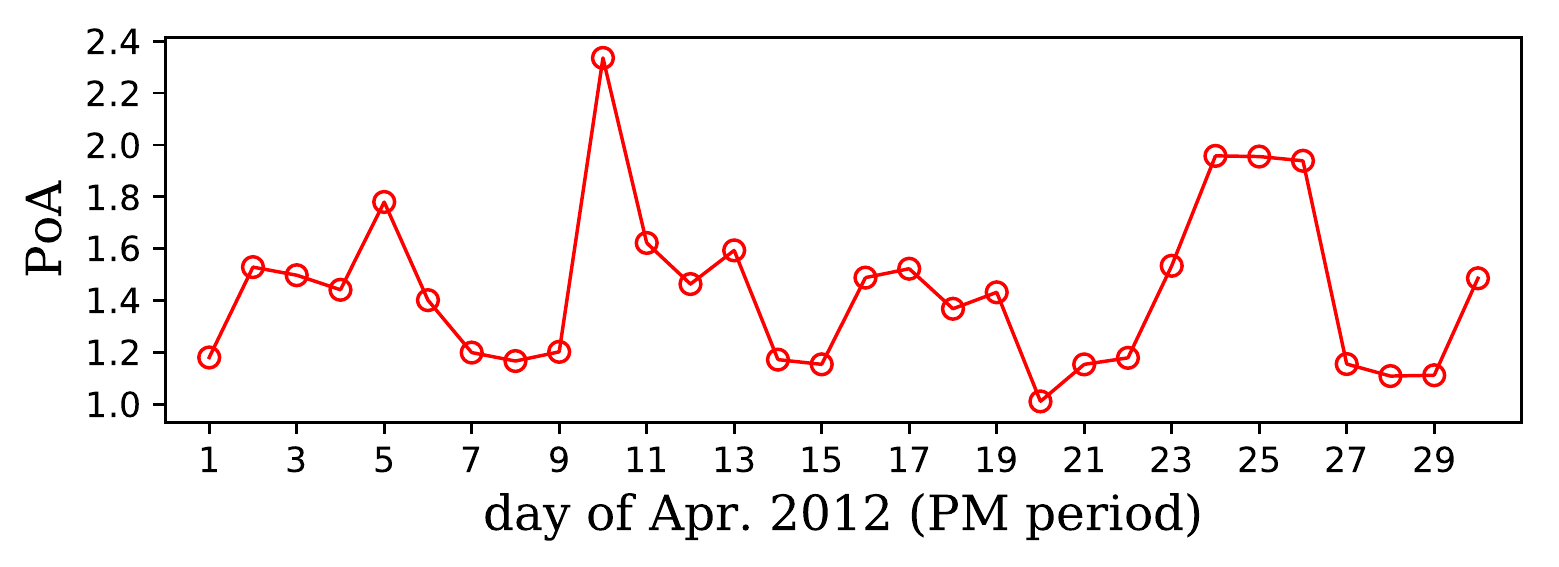}
		\caption{for $\scrI_2$}
		\label{poa-ifac}
	\end{subfigure} 
	\begin{subfigure}[b]{0.49\textwidth}
		\includegraphics[width=\textwidth]{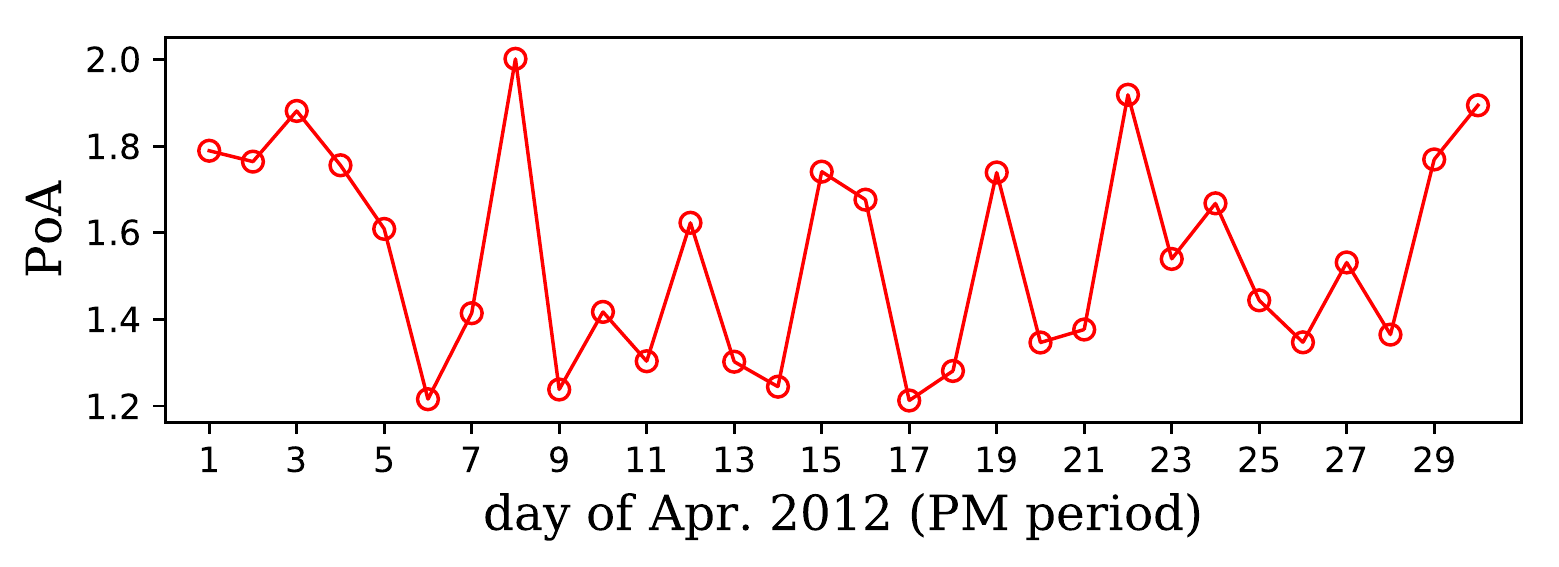}
		\caption{for $\scrI_3$}
		\label{poa-journal}
	\end{subfigure}   
	\caption{Daily PoAs for $\scrI_2$ and $\scrI_3$ (PM period for  Apr. 2012).}
	\label{poa-ieee}
\end{figure}

The PoA values for $\scrI_2$ shown in Fig. \ref{poa-ifac} have larger
variations than those for $\scrI_1$ in \cite{CDC16} and for $\scrI_3$
shown in Fig. \ref{poa-journal}; some are closer to 1 but some go beyond
2.2, meaning we have larger potential to improve the road network. {It
  is also seen that, although $\scrI_2$ is extracted from $\scrI_3$,
  there is no obvious correlation between the PoA values estimated for
  $\scrI_2$ and $\scrI_3$. To explain this, one should notice the fact
  that $\scrI_2$ is only a small subnetwork of $\scrI_3$, where the
  latter contains many more nodes/links/OD pairs (see
  Figs. \ref{sub-net-a-ext} and \ref{fig:zoneMA}). Specifically, in
  Fig.~\ref{fig:zoneMA} many more links have been added which
  significantly alter the feasible routing patterns relative to
  Fig.~\ref{sub-net-a-ext}. Thus, even though there may be correlations
  at the individual link flow level, once we add links and then
  aggregate over all links, any correlation is likely weakened or
  lost. Moreover, the social optimization problems solved to obtain the
  denominator of the PoA ratio in \eqref{poa-def} are very different
  since the subnetwork topologies are different.}  However, when taking
the average of the PoA values for all 30 days of Apr. 2012, all
$\scrI_1$, $\scrI_2$, and $\scrI_3$ result in an average PoA
approximately equal to 1.5, meaning we can gain an efficiency
improvement of about 50\%; thus, the results are consistent.

		\subsection{Results from sensitivity analysis}  \label{sec:sensResult}
		 
\begin{figure}
	\centering
		\includegraphics[width=0.6\textwidth]{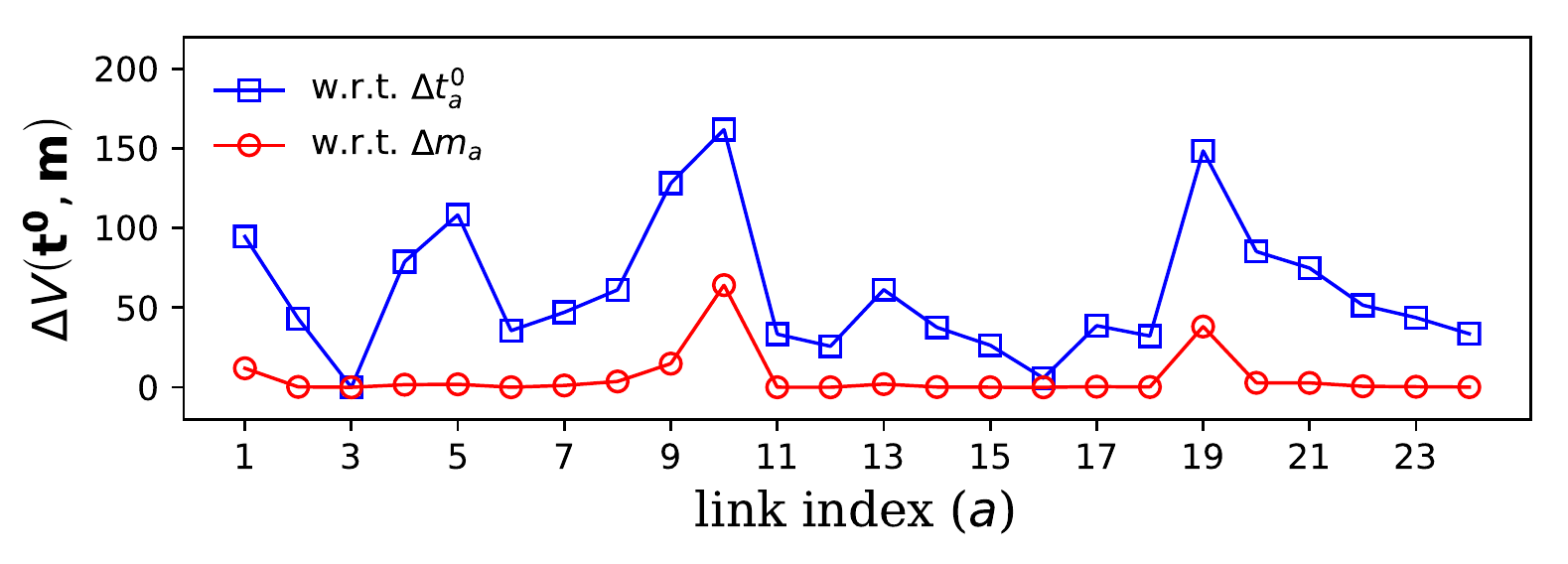}
	\caption{Sensitivity analysis (finite difference approximation) results for $\scrI_1$;  AM period of Apr. 2012.}
	\label{fig:sens-fda}
\end{figure}

		Investigating the AM peak period of Apr. 2012 for $\scrI_1$, instead of directly applying the formulae \eqref{cdc16-s1} and \eqref{cdc16-s2}, we calculate the two quantities defined in \eqref{sensi-freetime} and \eqref{sensi-capacity}, and plot the results in Fig. \ref{fig:sens-fda}, where the blue (resp., red) curve indicates the quantity $\Delta V\left( {{\bt^0},\bm;\Delta t_{a}^0} \right)$ (resp., $\Delta V\left( {{\bt^0},\bm;\Delta m_{a}} \right)$) for each and every link of $\scrI_1$. 
		It is seen from Fig. \ref{fig:sens-fda} that the largest four values of $\Delta V\left( {{\bt^0},\bm;\Delta t_{a}^0} \right)$ (resp., $\Delta V\left( {{\bt^0},\bm;\Delta m_{a}} \right)$) correspond to links {10, 19, 9, and 5 (resp., 10, 19, 9, and 1)}. This suggests that, during the AM peak
		period of Apr. 2012, the transportation management department
		could have most efficiently reduced the objective function value of the TAP \eqref{cdc16-tap}, thus mitigating congestion, by
		taking actions with priorities on these links (e.g., improving road
		conditions to reduce the free-flow travel time for links {10, 19, 9, and 5}, and increasing the number of lanes to enlarge the flow capacity
		for links {10, 19, 9, and 1}).

\subsection{Results from meta analysis}	

We conduct meta analysis for $\scrI_3$, under the user-centric routing
policy vs. the system-centric one. Our analysis includes the zone costs,
the maximum/minimum link flows, and the link-specific congestion.

\subsubsection{Meta analysis for zone costs}

Let $\scrA_{3}^i$ denote the set of links related to zone $i$ of
$\scrI_3$ (each link in $\scrA_{3}^i$ has at least one node contained in
zone $i$). Then, the total users' travel latency cost for zone $i$ is
defined as
$${C_i} = \sum\limits_{a \in {\scrA_{3}^i}} {{{x_a}} {{t_a}\left( x_a \right)} }. $$

\begin{figure}[h]  
	\centering
	\begin{subfigure}[b]{0.49\textwidth}
		\includegraphics[width=\textwidth]{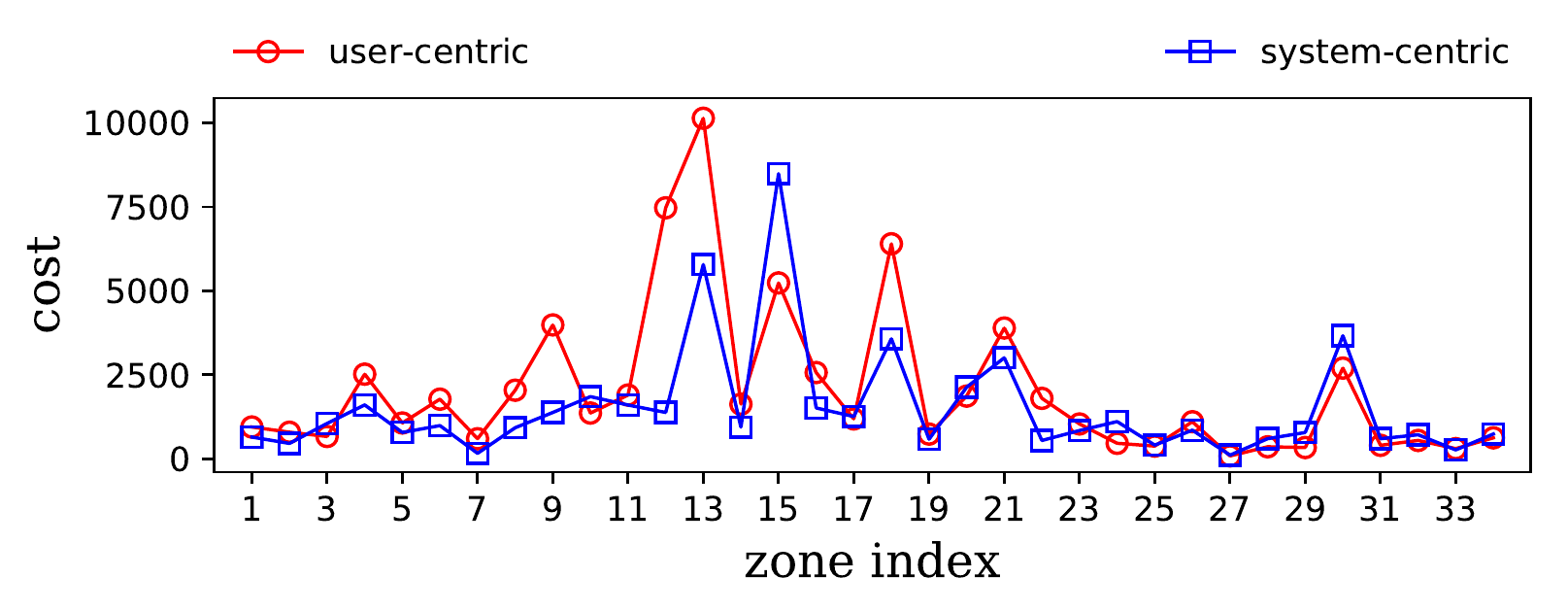}
		\caption{weekday (4/18/2012)}
		\label{zoneCostWeekday}
	\end{subfigure} 
	\begin{subfigure}[b]{0.49\textwidth}
		\includegraphics[width=\textwidth]{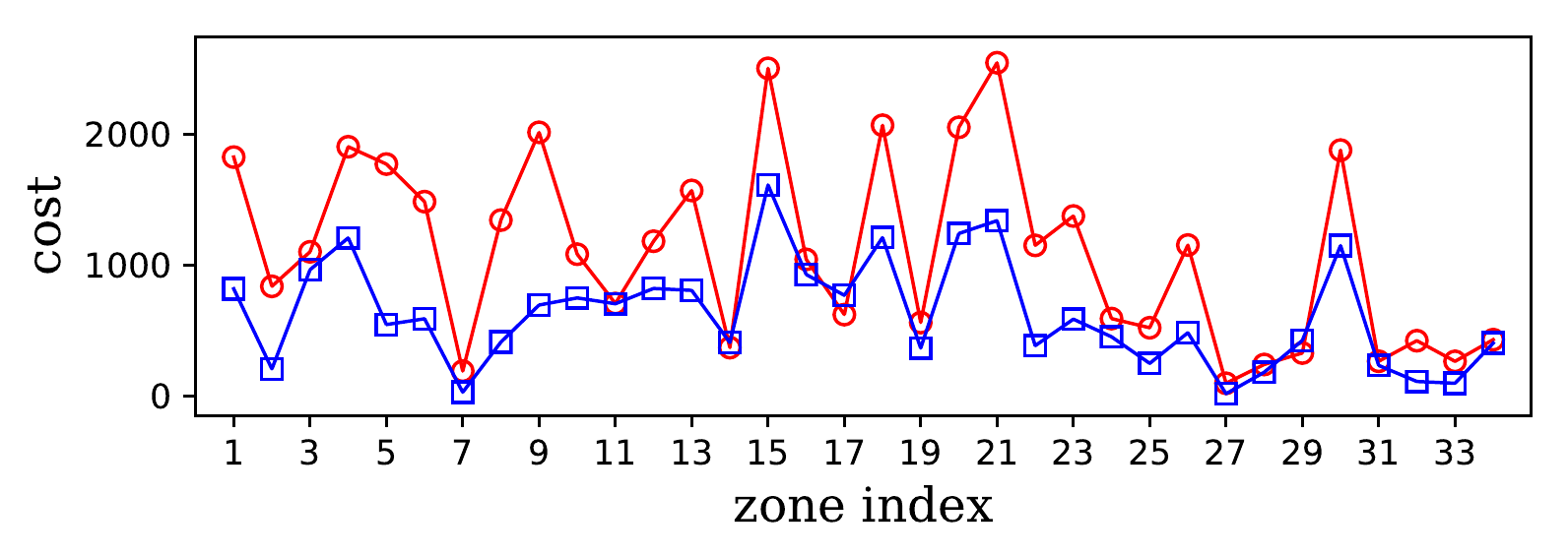}
		\caption{weekend (4/15/2012)}
		\label{zoneCostWeekend}
	\end{subfigure}   
	\caption{Zone costs under user-centric vs. system-centric routing policy (PM period of Apr. 2012).}
	\label{fig:zoneCost}
\end{figure}

We consider two scenarios, one corresponding to the PM peak period of a
typical weekday (Wednesday, 4/18/2012) and the other the PM period of a
typical weekend (Sunday, 4/15/2012). The zone costs under the
user-centric (resp., system-centric) routing policy are visualized in
Fig. \ref{zoneCostWeekday} (resp., \ref{zoneCostWeekend}).  Three
observations can be made: (\rmnum{1}) Overall, most zone costs would be
reduced when switching from the user-centric routing policy to the
system-centric one. (\rmnum{2}) In general, the zone costs for weekends
are less than their counterparts for weekdays; this is consistent with
intuition. (\rmnum{3}) The decrease seems more consistent for all zones
during weekends than during weekdays, suggesting it is easier to
optimize the network during weekends; this is again consistent with
intuition.

\subsubsection{Meta analysis for maximum/minimum link flows}

\begin{figure} 
	\centering
	\begin{subfigure}[b]{0.49\textwidth}
		\includegraphics[width=\textwidth]{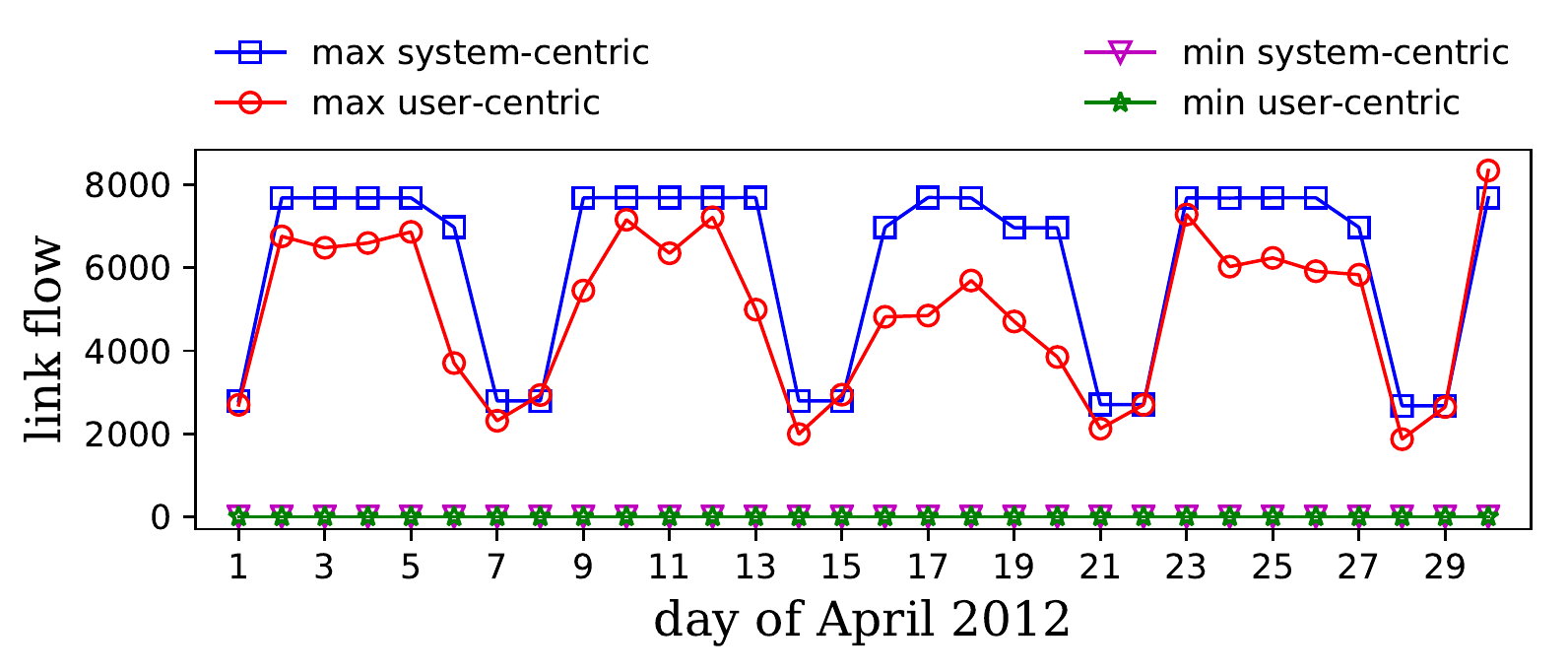}
		\caption{}
		\label{maxMinLinkFlow}
	\end{subfigure} 
	\begin{subfigure}[b]{0.49\textwidth}
		\includegraphics[width=\textwidth]{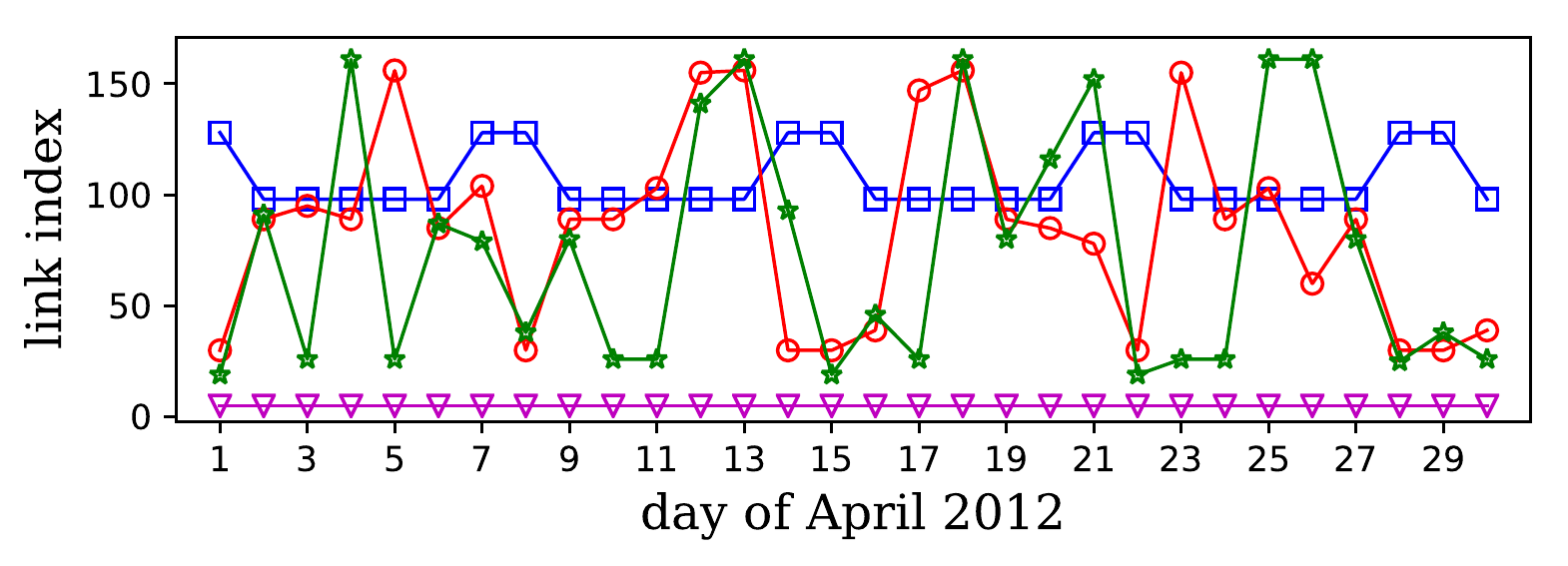}
		\caption{}
		\label{maxMinLinkFlowIdx}
	\end{subfigure}   
	\caption{Maximum/minimum link flows and the corresponding link indices under user-centric vs. system-centric routing policy (PM period of Apr. 2012).}
	\label{fig:maxMinLinkFlow}
\end{figure}

The maximum/minimum link flows for the PM peak period of each and every
day of Apr. 2012 are plotted in Fig. \ref{maxMinLinkFlow}, and the
corresponding link indices are shown in Fig. \ref{maxMinLinkFlowIdx}. A
major observation, based on Fig. \ref{maxMinLinkFlow}, is that the
maximum link flow values would increase for most of the days when
switching the routing policy from the user-centric one to the
system-centric one, which is desirable. {In addition, it is seen that, among the entire month (April 2012), both the maximum link flows under the two routing policies have a weekly periodic distribution; this is consistent with intuition.}

\subsubsection{Meta analysis for link congestion}

For any given link $a$, we define its \emph{Congestion Metric (CM)}
\cite{aftabuzzaman2007measuring} as the ratio of the $\textit{travel
  time}$ to $\textit{free-flow travel time}$: 
\begin{align}
{\text{CM}_a} \defeq \frac{{{t_a}\left( {{x_a}} \right)}}{{t_a^0}} = f\left( {\frac{{{x_a}}}{{{m_a}}}} \right),
\label{cm}
\end{align}
where $f(\cdot)$ is the cost function that we have estimated. By this definition, we always have $\text{CM}_a \geq 1$.

\begin{figure*}
	\centering
		\includegraphics[width=\textwidth]{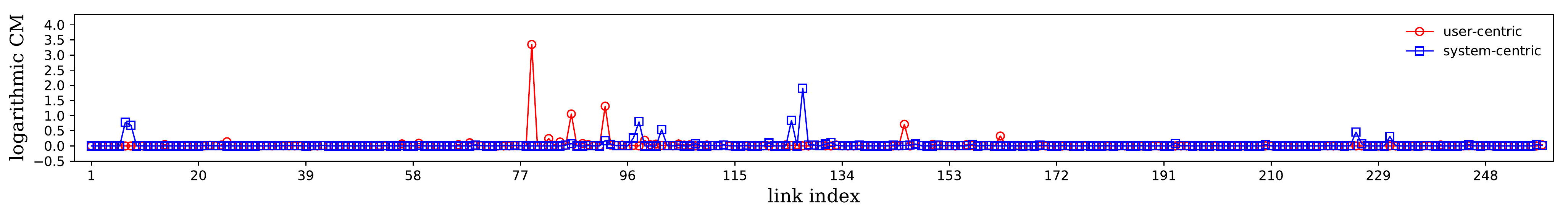}
	\caption{Link congestion under user-centric vs. system-centric routing policy (PM period of 4/18/2012).}
	\label{linkCongWeekday}
\end{figure*}

We first consider a PM peak period scenario for a typical workday
(Wednesday, 4/18/2012). The CM values of all the 258 links are plotted
in Fig. \ref{linkCongWeekday} {in a logarithmic scale (base 2)}. It
is seen that, for some links (indexed with 79, 92, and 86) the CM value
is \emph{significantly higher} ($\text{gap} > 1$) under the user-centric
routing policy than under the system-centric one. There are some links
for which we have the opposite, but, overall, the CM peak is reduced
under the system-centric policy. We then investigate a PM period
scenario for a typical weekend (Sunday, 4/15/2012), and find that all
the CM values for this scenario are very close to 1, meaning there was
almost no congestion for all links; we have omitted the weekend CM plot
for economy of space.

\section{Strategies for PoA reduction}  \label{sec:strat}

After quantifying the PoA, a natural question we must answer is the
following: How can we reduce the PoA for a given transportation network?
We propose three practical strategies for reducing the PoA, especially
when $\text{PoA} \gg 1$.

First, by taking advantage of the rapid emergence of Connected
Automated Vehicles (CAVs), it has become feasible to automate routing decisions, thus solving a \emph{system-centric forward
problem} (cf. \eqref{obj_soc}) in which all CAVs (bypassing driver decisions) cooperate to optimize the overall system performance. 

Second, we
propose a modification to existing GPS navigation algorithms recommending to all drivers socially optimal routes, which could be implemented by making use of \eqref{user-social}. In particular, we can solve the \emph{user-centric forward problem} \eqref{cdc16-tap}, embedded in a typical GPS navigation application, with $t_a(\cdot)$ replaced by $\overline t_a(\cdot)$, whose common cornerstone part, $f(\cdot)$, is estimated using \eqref{cdc16-costEstimator}. It is worth pointing out that some existing work simply took $f(\cdot)$ to be the Bureau of Public Roads (BPR)'s \cite{branston1976link} empirical polynomial function $f(z) = 1 + 0.15z^4,\,\forall z \ge 0$, which would not be as accurate.

Finally, our
sensitivity analysis results provide the means to prioritize road segments for specific interventions that can mitigate congestion.

\section{Conclusions and Future Work}  \label{sec:conc}

In this paper, we assess the efficiency of transportation networks under a selfish user-centric routing policy as opposed to a socially-optimal system-centric routing policy. To that end, we define and quantify the Price of Anarchy (PoA) and propose possible strategies to reduce it. All the procedures involved are data-driven, thus having the capability of dynamically optimizing any given transportation network (by using the data collected in real-time manner), in terms of reducing the PoA (especially when $\text{PoA} \gg 1$) such that it gets as close to 1 as possible.

{We must keep in mind that, due to unavoidable inaccuracies in data
  and modeling, all the numerical results shown in Sec. \ref{Sec:Rsults}
  are only estimates. In particular, the speed-to-flow conversion model
  that we use (Greenshield's model) is a macroscopic model with
  naturally limited accuracy, the GLS method that we leverage also is
  based on an approximation, and the MSA subroutine in Alg. 1 is an
  approximate scheme.}

{In terms of the computational challenges of our proposed
  approaches, we encountered numerical difficulties when solving
  \eqref{qp2} and \eqref{qp3} to obtain OD demands for large-sized (say
  a network like $\scrI_3$) networks. However, we subsequently developed
  a simplification procedure by considering only the
  {\emph{fastest}} route for each OD pair, thus successfully
  resolving this issue. We conducted case studies on a workstation with
  24 GB memory and a 12-core Intel Core i5 CPU, and for the largest
  network ($\scrI_3$) that we investigated, the total CPU time
  (including estimating OD demands, recovering link latency cost
  functions, adjusting OD demands, solving for socially-optimal flows,
  and finally calculating PoA values) is about 10 hours. The total CPU
  times for $\scrI_1$ and $\scrI_2$ are about 30 minutes and 2 hours,
  respectively. We note that the most time-consuming task is adjusting
  OD demands using Alg. \ref{alg:demandAdjustment}. However, it is seen
  that Steps 2-4 of Alg. \ref{alg:demandAdjustment} and the MSA
  subroutine can easily benefit from parallel computing. Thus,
  scalability can be further improved through parallel
  computation. Moreover, following an approach of ``divide and
  conquer,'' several decomposition methods could possibly also be
  leveraged as we move to larger networks; the difficulty lies in how to
  reasonably ``merge'' results derived for subnetworks so as to obtain
  the final result for the whole network.}

Our ongoing work includes extending the PoA analysis and reduction
framework from single-class to multi-class transportation networks. We
have recently obtained results for the \emph{multi-class user-centric
  inverse problem} \cite{multi-class-inverse-equilibrium-cdc2017}, which
paves the way for data-driven PoA estimation in these networks. We are
also considering alternative models/methods to improve the accuracy in
the PoA evaluation. In addition, it is of interest to consider jointly
estimating/adjusting the OD demand matrices and recovering the travel
latency cost functions.

\section*{Acknowledgments}

The authors would like to thank the Boston Region Metropolitan Planning Organization, and Scott Peterson in particular, for supplying the EMA data and providing us invaluable clarifications throughout our work. 

\bibliographystyle{IEEEtran}





\end{document}